\documentclass[11pt,fleqn,twoside]{article}
\usepackage{amsfonts,amssymb,latexsym}
\makeatletter
\newcommand{\prava}[1]{\small\it
\begin{flushleft}
Copyright \copyright \ 1999 by  #1
\end{flushleft}}

\newcommand{\name}[1]{\begin{flushleft}
                       \LARGE \bf #1
                       \end{flushleft}\vspace{-3mm}}

\newcommand{\Author}[1]{\begin{flushleft}
                       \it #1 \end{flushleft}}

\newcommand{\Adress}[1]{\begin{flushleft}
                       \it #1 \end{flushleft}}

\newcommand{\Date}[1]{\begin{flushleft}
                      \small  \it #1 \end{flushleft}}

\newcommand{\ehkol}{Author \ name}
\newcommand{\ohkol}{Article \ name}
\renewcommand{\@evenhead}{
\hspace*{-3pt}\raisebox{-15pt}[\headheight][0pt]{\vbox{\hbox to \textwidth 
{\thepage \hfil \ehkol}\vskip4pt \hrule}}}
\renewcommand{\@oddhead}{
\hspace*{-3pt}\raisebox{-15pt}[\headheight][0pt]{\vbox{\hbox to \textwidth 
{\ohkol \hfil \thepage}\vskip4pt\hrule}}}
\renewcommand{\@evenfoot}{}
\renewcommand{\@oddfoot}{}

\newcommand{\be}{\begin{equation}}
\newcommand{\ee}{\end{equation}}
\newcommand{\ba}{\hspace*{-5pt}\begin{array}}
\newcommand{\ea}{\end{array}}

\newcommand{\ds}{\displaystyle}
\makeatother

\def\tfr#1#2{{\textstyle{#1\over#2}}}

\def\d{{\rm d}}
\def\e{{\rm e}}
\def\i{{\rm i}}

\def\~#1{{\bf\widetilde{\mit#1}}}
\def\^#1{{\bf\hat{\mit#1}}}
\def\=#1{{\bf\bar{\mit#1}}}

\def\ch{Camassa-Holm}

\def\Rl{{\rm I\kern-.25em R}}  
\def\Cx{{\rm C                 
        \kern-.49em\vrule depth-0.06ex width0.05em height1.44ex
        \kern-.15em\vrule depth-0.35ex width0.10em height1.22ex
        \kern .55em}}

\def\acs{arbitrary constants}
\def\isc{invariant surface condition}

\def\af{arbitrary function}
\def\afs{arbitrary functions}
\def\inl{inf\/initesimal}
\def\inls{inf\/initesimals}
\def\wrt{with respect to}

\def\deqs{determining equations}
\def\pde{partial dif\/ferential equation}
\def\pdes{partial dif\/ferential equations}

\def\wz{\frac{\d w}{\d z}}
\def\wzz{\frac{\d^2 w}{\d z^2}}
\def\wzzz{\frac{\d^3 w}{\d z^3}}
\def\wzzzz{\frac{\d^4 w}{\d z^4}}
\def\Wz{\frac{\d W}{\d z}}

\def\wx{\frac{\d w}{\d x}}
\def\wxx{\frac{\d^2 w}{\d x^2}}
\def\wxxx{\frac{\d^3 w}{\d x^3}}
\def\wxxxx{\frac{\d^4 w}{\d x^4}}

\def\wtt{\frac{\d^2 w}{\d t^2}}
\def\wwtt{\frac{\d^2 \~{w}}{\d t^2}}

\def\gt{\frac{\d g}{\d t}}
\def\gtt{\frac{\d^2 g}{\d t^2}}
\def\gttt{\frac{\d^3 g}{\d t^3}}

\def\Htt#1{\frac{\d^2 H_{#1}}{\d t^2}}

\def\beq{\begin{equation}}
\def\eeq{\end{equation}}
\def\bear{\begin{eqnarray}}
\def\ear{\end{eqnarray}}
\def\bearn{\begin{eqnarray*}}
\def\earn{\end{eqnarray*}}

\def\sen{2} 
\newcount\redno
\def\redn#1{\global\advance\redno by 1
           \Reduction{\sen.\the\redno}
           \expandafter \xdef\csname #1\endcsname {\sen.\the\redno}\relax }

\def\Reduction#1{\vspace{0.2cm}
\noindent{\bf Reduction\ #1}}

\def\eqtab#1#2{(#1#2)}

\def\pac{Clarkson P.A.}
\def\liz{Mansf\/ield E.L.}
\def\hs{Segur H.}
\def\mja{Ablowitz M.J.}

\begin{document}
\setcounter{page}{66}
\thispagestyle{empty}

\renewcommand{\ehkol}{P.A.~Clarkson and T.J.~Priestley}
\renewcommand{\ohkol}{Symmetries of a Class of Nonlinear Fourth Order Partial
Dif\/ferential Equations}

\begin{flushleft}
\footnotesize \sf
Journal of Nonlinear Mathematical Physics \qquad 1999, V.6, N~1,
\pageref{clarkson-fp}--\pageref{clarkson-lp}.
\hfill {\sc Article}
\end{flushleft}

\vspace{-5mm}

\renewcommand{\footnoterule}{}
{\renewcommand{\thefootnote}{} \footnote{\prava{P.A.~Clarkson and
T.J.~Priestley}}} 

\name{Symmetries of a Class of Nonlinear Fourth Order Partial
Dif\/ferential Equations} \label{clarkson-fp}

\Author{Peter A.~CLARKSON and Thomas J.~PRIESTLEY}

\Adress{Institute of Mathematics and Statistics, University of Kent
at Canterbury,\\
 Canterbury, CT2 7NF, UK}

\Date{Received September 01, 1998}

\begin{abstract}
\noindent
In this paper we study symmetry reductions of a class of nonlinear
fourth order \pdes 
\be
u_{tt} = \left(\kappa u + \gamma u^2\right)_{xx} + u u_{xxxx}
+\mu u_{xxtt}+\alpha u_x u_{xxx} + \beta u_{xx}^2, 
\ee
where $\alpha$, $\beta$, $\gamma$, $\kappa$ and $\mu$ are \acs.
This equation may be thought of as a fourth order analogue of a
generalization of the \ch\ equation, about which there has been
considerable recent interest. Further equation (1) is a
``Boussinesq-type'' equation which arises as a model of 
vibrations of an anharmonic mass-spring chain and admits both
``compacton'' and conventional solitons.
 A catalogue of symmetry reductions for equation (1) is obtained using
the classical Lie method and the nonclassical method due to Bluman
and Cole. In particular we obtain several reductions using the
nonclassical method which are {\it not} obtainable 
through the classical method.
\end{abstract}

\setcounter{equation}{0}

\section{Introduction}
In this paper we are concerned with symmetry reductions
of the nonlinear fourth order \pde\ given by
\beq\Delta \equiv u_{tt} - \left(\kappa u + \gamma u^2\right)_{xx} -
u u_{xxxx} - \mu u_{xxtt}- \alpha u_x u_{xxx} - \beta u_{xx}^2 =0,
\label{fulleqn}
\eeq 
where $\alpha$, $\beta$, $\gamma$, $\kappa$ and $\mu$ are \acs. 
This equation may be thought of as an alternative to a
generalized \ch\ equation (cf.\ \cite{refCMP} and
the references therein)
\beq
u_t -\epsilon u_{xxt} +2\kappa u_x= u u_{xxx}
+\alpha u u_x +\beta u_x u_{xx}. \label{gch}
\eeq
This is analogous to the Boussinesq equation \cite{refBi,refBii}
\beq
u_{tt}  = \left(u_{xx} + \tfr12 u^2\right)_{xx} \label{bq}
\eeq
which is a soliton equation solvable by inverse scattering
\cite{refAH,refCi,refCii,refDTT,refZ},
being an alternative to the Korteweg-de Vries (KdV) equation
\beq 
u_t = u_{xxx} + 6uu_x \label{kdv}
\eeq
another soliton equation, the f\/irst to be solved by inverse
scattering~\cite{refGGKM}. 

Two special cases of (\ref{fulleqn}) have appeared recently in the
literature both of which model the motion of a dense chain
\cite{refRosa}. The f\/irst is obtainable via the transformation 
\[
(u,x,t)\mapsto (2\varepsilon \alpha_3 u+\varepsilon \alpha_2, x,t)
\]
with the appropriate change of parameters, 
to yield 
\beq
u_{tt} = \left(\alpha_2 u+ \alpha_3 u^2\right)_{xx} + \varepsilon \alpha_2
 u_{xxxx} + 2 \varepsilon \alpha_3 \left[ u u_{xxxx} +
 2 u_{xx}^2 +3 u_{x} u_{xxx} \right] \label{rone}
\eeq
with $\varepsilon>0$. This equation can be thought of as the
Boussinesq equation (\ref{bq}) appended with a nonlinear dispersion.
It admits both conventional solitons and compact solitons often
called ``compactons''. Compactons are solitary waves with a compact
support (cf.\ \cite{refRosa,refRosb,refRosc,refRH}). The compact
structures take the form 
\beq
u(x,t) = 
\left\{
\begin{array}{ll}
\ds\frac{3 c^2 - 2\alpha_2}{2 \alpha_3} \cos^2
\left\{ (12 \varepsilon)^{-1/2} (x-ct)\right\}, &
\mbox{if} \quad |x-ct|\le 2\pi, 
\vspace{2mm}\\
 0, & \mbox{if} \quad |x-ct|> 2\pi.
\end{array}\right.
\label{componei}
\eeq
or
\beq 
u(x,t) = \left\{
\begin{array}{ll}
A \cos \left\{ (3 \varepsilon)^{-1/2} \left[x-
 \left( \tfr23 \alpha_3 \right)^{1/2} t \right]\right\} ,&
\mbox{if} \quad |x-ct|\le 2\pi, 
\vspace{2mm}\\
 0, & \mbox{if} \quad |x-ct|> 2\pi.
\end{array}\right.
\label{componeii}
\eeq
These are ``weak'' solutions as they do not possess the necessary
smoothness at the edges, however this would appear not to af\/fect the
robustness of a compacton \cite{refRosa}. Numerical experiments seem
to show that compactons interact elastically, reemerging with exactly
the same coherent shape \cite{refRH}. See \cite{refLOR} for a recent
study of non-analytic solutions, in particular compacton solutions,
of nonlinear wave equations. 

The second equation is obtained from the scaling transformation 
\[
(u,x,t) \mapsto \left(2 \alpha_3 u/\varepsilon,\sqrt \varepsilon\,
x,t\right), 
\]
again with appropriate parameterisation,
\beq 
u_{tt} = \left(\alpha_2 u+ \alpha_3 u^2\right)_{xx} + \varepsilon
u_{xxtt} + 2 \varepsilon \alpha_3 \left[ u u_{xxxx} +
 2 u_{xx}^2 +3 u_{x} u_{xxx} \right]
\label{rtwo}
\eeq
with $\varepsilon>0$. This equation, unlike
(\ref{rone}) is well posed. It also admits conventional solitons
and allows compactons like
\beq u(x,t) = 
\left\{
\begin{array}{ll}
\ds\frac{4 c^2 - 3\alpha_2}{2 \alpha_3} \cos^2
\left\{ (12 \varepsilon)^{-1/2} (x-ct)\right\}, &
\mbox{if} \quad |x-ct|\le 2\pi, 
\vspace{2mm}\\
0, & \mbox{if} \quad |x-ct|> 2\pi, 
\end{array}\right.
\label{comptwoi}
\eeq
or
\beq 
u(x,t) = \left\{
\begin{array}{ll}
A \cos \left\{ (3 \varepsilon)^{-1/2} \left[x-
\left( \tfr32 \alpha_2 \right)^{1/2} t \right]\right\} ,&
\mbox{if} \quad |x-ct|\le 2\pi, 
\vspace{2mm}\\
0, & \mbox{if} \quad  |x-ct|> 2\pi.
\end{array}\right. 
\label{comptwoii}
\eeq
These again are weak solutions, and are very similar to the previous
solutions: both (\ref{componeii}) and (\ref{comptwoii}) are solutions
with a variable speed linked to the amplitude of the wave, whereas
both (\ref{componei}) and (\ref{comptwoi}) are solutions with
arbitrary amplitudes, whilst the wave speed is f\/ixed by the
parameters of the equation.

The Fuchssteiner-Fokas-Camassa-Holm (FFCH) equation
\beq 
u_t -u_{xxt} +2\kappa u_x= u u_{xxx}
-3u u_x + 2u_x u_{xx}, \label{ffch}
\eeq
f\/irst arose in the work of Fuchssteiner and Fokas
\cite{refFucha,refFF} using a bi-Hamiltonian approach;
we remark that it is only implicitly written in \cite{refFF} --  see
equations (26e) and (30) in this paper -- though is explicitly
written down in \cite{refFucha}. It has recently been rederived by
Camassa and Holm \cite{refCH} from physical considerations as a model
for dispersive shallow water waves. In the case $\kappa=0$, it admits
an unusual solitary wave solution 
\[
u(x,t) = A\exp\left(-\vert x-ct \vert\right), 
\]
where $A$ and $c$ are arbitrary constants, which is called a
``peakon''. A Lax-pair \cite{refCH} and bi-Hamiltonian structure
\cite{refFF} have been found for the FFCH equation (\ref{ffch}) and
so it appears to be completely integrable. Recently the FFCH equation
(\ref{ffch}) has attracted considerable attention. In 
addition to the aforementioned, other studies include 
\cite{refCHH,refCona,refConb,refConc,refCS,refFokb,refFOR,refFuchb,refGP,refHone,refOR,refSchiff}.

Symmetry reductions and exact solutions have several dif\/ferent
important applications in the context of dif\/ferential equations.
Since solutions of \pdes\ asymptotically 
tend to solutions of lower-dimensional equations obtained by symmetry
reduction, some of these special solutions will illustrate important
physical phenomena. In particular, exact solutions arising from
symmetry methods can often be used ef\/fectively to study 
properties such as asymptotics and ``blow-up'' (cf.\ \cite{refG,refGDEKS}).
Furthermore, explicit solutions (such as those found by symmetry
methods) can play an important role in the design and testing of
numerical integrators; these solutions provide an important practical
check on the accuracy and reliability of such integrators (cf.\
\cite{refA,refS}). 

The classical method for f\/inding symmetry reductions of partial
dif\/ferential equations is the Lie group method of inf\/initesimal
transformations, which in practice is a two-step procedure (see \S~2
for more details). The f\/irst step is entirely algorithmic, 
though often both tedious and virtually unmanageable manually. As a
result, symbolic manipulation (SM) programs have been developed to
aid the calculations; an excellent survey of the dif\/ferent packages
available and a description of their strengths and 
applications is given by Hereman \cite{refH}. In this paper we use
the {\rm MACSYMA} package {\tt symmgrp.max} \cite{refCHW} to
calculate the determining equations. 
The second step involves heuristic integration procedures which have been
implemented in some SM programs and are largely successful, though
not infallible. Commonly, the overdetermined systems to be solved are
simple, and heuristic integration is both fast and ef\/fective.
However, there are occasions where heuristics 
can break down (cf.\ \cite{refMCi} for further details and examples).
Of particular importance to this study, is if the classical method is
applied to a \pde\ which contains arbitrary parameters, such as
(\ref{fulleqn}) or more generally, \afs. 
Heuristics usually yield the general solution yet miss those special
cases of the parameters and \afs\ where additional symmetries lie. In
contrast the method of dif\/ferential Gr\"obner bases (DGBs), which we
describe below, has proved ef\/fective in coping with such dif\/f\/iculties
(cf.\ \cite{refCMi,refCMP,refMCi,refMCii}). 

In recent years the nonclassical method due to Bluman and Cole
\cite{refBC} (in the sequel referred to as the ``nonclassical
method''), sometimes referred to as the ``method of partial
symmetries of the f\/irst type'' \cite{refV}, or the ``method of 
conditional symmetries'' \cite{refLW}, and the direct method due to
Clarkson and Kruskal \cite{refCK} have been used, with much success,
to generate many new symmetry reductions and exact solutions for
several physically signif\/icant \pdes\ that are not 
obtainable using the classical Lie method (cf.\ \cite{refPAC,refCLP}
and the references therein). The nonclassical method is a
generalization of the classical Lie method, whereas the direct method
is an ansatz-based approach which involves no group 
theoretic techniques. Nucci and Clarkson \cite{refNC} showed that for the
Fitzhugh-Nagumo equation the nonclassical method is more general than
the direct method, since they demonstrated the existence of a
solution of the Fitzhugh-Nagumo equation, obtainable using the
nonclassical method but not using the direct method. 
Subsequently Olver \cite{refOlverb} (see also \cite{refABH,refPucci})
has proved the general result that for a scalar equation, every
reduction obtainable using the direct method is also obtainable using
the nonclassical method. Consequently we use the nonclassical method
in this paper rather than the direct method. 

The method used to f\/ind solutions of the determining equations in
both the classical and nonclassical method is that of DGBs, def\/ined
to be a basis \ss\ of the dif\/ferential ideal generated by the system
such that every member of the ideal pseudo-reduces to zero with
respect to the basis \ss. This method provides a systematic framework
for f\/inding integrability and compatibility conditions of an 
overdetermined system of \pdes. It avoids the problems of inf\/inite
loops in reduction processes and yields, as far as is currently
possible, a ``triangulation'' of the system from which the solution
set can be derived more easily \cite{refCMi,refMF,refR,refRi}. In a
sense, a DGB provides the maximum amount of information possible
using elementary dif\/ferential and algebraic processes in f\/inite 
time.

In pseudo-reduction, one must, if necessary, multiply the expression
being reduced by dif\/ferential (non-constant) coef\/f\/icients of the
highest derivative terms of the reducing equation, so that the
algorithms used will terminate \cite{refMF}. In practice, such
coef\/f\/icients are assumed to be non-zero, and one needs to deal with
the possibility of them being zero separately. These are called
singular cases. 

The triangulations of the systems of \deqs\ for \inls\ arising in the
nonclassical method in this paper were all performed using the {\rm
MAPLE} package {\tt diffgrob2} \cite{refM}. This package was written
specif\/ically to handle nonlinear equations of polynomial type. All
calculations are strictly `polynomial', that is, there is no
division. Implemented there are the Kolchin-Ritt algorithm using
pseudo-reduction instead of reduction, and extra algorithms needed to
calculate a DGB (as far as possible using the current theory), for
those cases where the Kolchin-Ritt algorithm is not suf\/f\/icient
\cite{refMF}. The package was designed to be used interactively as 
well as algorithmically, and much use is made of this fact here. It
has proved useful for solving many fully nonlinear systems (cf.\
\cite{refCMi,refCMii,refCMiii,refCMiv,refCMP}). 

In the following sections we shall consider the cases $\mu=0$ and
$\mu\not=0$, when we set $\mu=1$ without loss of generality,
separately because the presence or lack of the 
corresponding fourth order term is signif\/icant. In \S~2 we f\/ind the
classical Lie group of symmetries and associated reductions of
(\ref{fulleqn}). In \S~3 we discuss 
the nonclassical symmetries and reductions of (\ref{fulleqn}) in the
generic case. In \S~4 we consider special cases of the the
nonclassical method in the so-called $\tau=0$ case; in full
generality this case is somewhat intractable. In \S~5 we discuss our results.

\section{Classical Symmetries}
To apply the classical method we consider the one-parameter Lie group of inf\/initesimal 
trans\-formations in ($x,t,u$) given by
\be \label{trans}
\ba{l}
x^* = x+\varepsilon \xi (x,t,u) + {\cal O}\left(\varepsilon^2\right), 
\vspace{2mm}\\
t^* = t+\varepsilon \tau (x,t,u) + {\cal
O}\left(\varepsilon^2\right),
\vspace{2mm}\\ 
u^* = u+\varepsilon \phi (x,t,u) + {\cal
O}\left(\varepsilon^2\right), 
\ea
\ee
where $\varepsilon$ is the group parameter. Then one requires that
this transformation leaves invariant the set
\beq 
S_{\Delta} \equiv \{ u(x,t) : \Delta =0 \} \label{sdel}
\eeq
of solutions of (\ref{fulleqn}). This yields an overdetermined,
linear system of equations for the inf\/initesimals $\xi (x,t,u)$,
$\tau (x,t,u)$ and $\phi (x,t,u)$. The associated Lie algebra is
realised by vector f\/ields of the form 
\beq
{\bf v} = \xi (x,t,u) \frac{\partial}{\partial x} +\tau (x,t,u)
\frac{\partial}{\partial t} + \phi (x,t,u) \frac{\partial}{\partial
u}. \label{vecfi}
\eeq 
Having determined the inf\/initesimals, the symmetry variables are
found by solving the characteristic equation
\beq
{\d x \over \xi (x,t,u)} = {\d t \over \tau (x,t,u)}=
{\d u \over \phi (x,t,u)}, \label{chareq}
\eeq
which is equivalent to solving the \isc\
\beq
\psi \equiv \xi (x,t,u) u_x+\tau (x,t,u) u_t -
\phi (x,t,u) =0. \label{invsc}
\eeq
The set $S_{\Delta}$ is invariant under the transformation
(\ref{trans}) provided that 
\[
 {\rm pr}^{(4)} {\bf v} (\Delta)|_{\Delta \equiv 0} =0,
\]
where ${\rm pr}^{(4)}{\bf v}$ is the fourth prolongation of the
vector f\/ield (\ref{vecfi}), which is given explicitly in terms of
$\xi$, $\tau$ and $\phi$ (cf.\ \cite{refO}). 
This procedure  yields the determining equations. There are two
cases to consider, (i) $\mu=0$ and (ii) $\mu\not=0$.

\subsection{Case (i) $\mu=0$.}

In this case we generate 15 \deqs, using the {\rm MACSYMA} package
\linebreak {\tt symmgrp.max}. 
\[
\ba{l}
\tau_{u}=0,\qquad \tau_{x}=0, \qquad \xi_{u}=0, 
\qquad \phi_{uu}=0, \qquad \xi_{t}=0,
\vspace{1mm}\\
 \alpha (\phi_{u} u - \phi)=0, \qquad \beta (\phi_{u} u - \phi)=0, 
 \qquad 2 \phi_{tu} - \tau_{tt}=0, 
\vspace{1mm}\\
4 \phi_{xu} u - 6 \xi_{xx} u + \alpha \phi_{x}, \qquad
2 \tau_{t} u - 4 \xi_{x} u + \phi=0,
\vspace{1mm}\\
4 \beta \phi_{xu} + 3 \alpha \phi_{xu} - 2 \beta \xi_{xx} 
- 3 \alpha \xi_{xx}=0,
\vspace{1mm}\\
\phi_{tt} - \phi_{xxxx} u - 2 \gamma \phi_{xx} u - \kappa
\phi_{xx}=0,
\vspace{1mm}\\
3 \alpha \phi_{xxu} u + 2 \gamma \phi_{u} u + 4 \xi_{x} \gamma u 
- \alpha \xi_{xxx} u- 2 \gamma \phi=0,
\vspace{1mm}\\
6 \phi_{xxu} u^{2} + 4 \xi_{x} \gamma u^{2} - 4 \xi_{xxx} u^{2}
 + 2 \beta \phi_{xx} u + 2 \xi_{x} \kappa u - \kappa \phi=0,
\vspace{1mm}\\
4 \phi_{xxxu} u + 4 \gamma \phi_{xu} u - 2 \xi_{xx} \gamma u - \xi_{xxxx} u
+ \alpha \phi_{xxx} + 4 \gamma \phi_{x} + 2 \kappa \phi_{xu} -
\xi_{xx} \kappa=0, 
\ea
\]
and then use {\tt reduceall} in {\tt diffgrob2} to simplify
them to the following system
\[
\ba{l}
\xi_{u}=0, \qquad\xi_{t}=0, \qquad \gamma (14 \beta + 9 \alpha)
\xi_{x}=0, \qquad 
\tau_{u}=0, \vspace{1mm}\\ 
\gamma \kappa (14 \beta + 9 \alpha) \tau_{t}=0, \qquad \tau_{x}=0,\qquad
\gamma \kappa (14 \beta + 9 \alpha) \phi =0.
\ea
\]
Thus we have special cases when $\gamma=0$, $\kappa=0$ and/or $14
\beta + 9 \alpha=0$. The latter condition provides nothing dif\/ferent
unless we specialize further and consider the special case when
$\alpha=-\tfr52$ and $ \beta=\tfr{45}{28}$. We continue to use {\tt
reduceall} in {\tt diffgrob2} for the various combinations and it
transpires that there are only four combinations which yield
dif\/ferent \inls. Where a parameter is not included it is presumed to
be arbitrary. 
\[
\ba{lll} 
\mbox{(a)} & \kappa=0, & \xi_u=0,\quad \xi_{t}=0,
\quad\xi_{x}=0,\quad \tau_{u}=0, \vspace{2mm}\\
& &\tau_{tt}=0,\quad\tau_{x}=0,\quad 2\tau_t u+\phi=0. 
\vspace{2mm}\\
\mbox{(b)} & \gamma=0, & \xi_u=0,\quad \xi_{t}=0, \quad\xi_{xx}=0,
\quad \tau_{u}=0, \vspace{2mm}\\
& &\xi_x-\tau_t=0,\quad \tau_x=0,\quad 2\xi_x u+\phi=0. 
\vspace{2mm}\\
\mbox{(c)} & \gamma=\kappa=0, & \xi_u=0,\quad \xi_{t}=0, \quad
\xi_{xx}=0,\quad \tau_{u}=0,
\vspace{2mm}\\
& &\tau_{tt}=0, \quad \tau_{x}=0,\quad 2\tau_t u-4\xi_x u-\phi=0. 
\vspace{2mm}\\
\mbox{(d)} & \alpha=-\tfr52,\quad
\beta=\tfr{45}{28},\quad\gamma=\kappa=0, \quad
& \xi_u=0,\quad \xi_{t}=0, \quad\xi_{xxx}=0,\quad \tau_{u}=0,
\vspace{2mm}\\
& &\tau_{tt}=0, \quad \tau_{x}=0,\quad 2\tau_t u-4\xi_x u-\phi=0. 
\ea
\]

\noindent
Hence we obtain the following \inls.

\def\icinl{2.1}
\begin{center}
{\bf Table \icinl}

\medskip

\begin{tabular}{|ccccr|}
\hline & & & & \\[-5pt]
Parameters & $\xi$ & $\tau$ & $\phi$ & 
\\[5pt] \hline & & & & 
\\[-5pt]
& $c_1$ & $c_2$ & $0$ & (\icinl{i}) \\[5pt] 
$\kappa=0$ & $c_1$ & $c_3 t +c_2$ & $-2 c_3 u$ & (\icinl{ii}) \\[5pt] 
$\gamma=0$ & $c_3 x+c_1$ & $c_3 t +c_2$ & $2 c_3 u$ & (\icinl{iii}) \\[5pt] 
$\gamma=\kappa=0$ & $c_4 x+c_1$ & $c_3 t +c_2$ & $(4 c_4 - 2c_3) u$ &  
(\icinl{iv}) \\[5pt]
$\ba{c} \alpha=-\tfr52,\quad\beta=\tfr{45}{28}\\[5pt]
\gamma=\kappa=0 \ea $ &
$c_5 x^2+c_4 x+c_1$ & $c_3 t +c_2$ & $[4(2c_5 x+c_4)-2c_3] u$ &
(\icinl{v})\\[5pt] \hline 
\end{tabular}
\end{center}

\noindent
where $c_1,c_2,\ldots,c_5$ are \acs.

Solving the \isc\ yields the following seven dif\/ferent canonical reductions: 

\redn{ici} $\alpha$, $\beta$, $\gamma$ and $\kappa$ arbitrary. If in
(\icinl{i}--\icinl{v}) $c_3=c_4=c_5=0$, then we may set $c_2=1$ without
loss of generality. Thus we obtain the reduction
\[
u(x,t) = w(z), \qquad z=x-c_1 t
\]
where $w(z)$ satisf\/ies
\[
\left( \kappa -c_1^2\right)
\wzz + 2 \gamma \left[w \wzz + \left(\wz\right)^2\right]+
w \wzzzz + \alpha \wz \wzzz + \beta \left(\wzz\right)^2 =0. 
\]

\redn{icii} $\alpha$, $\beta$ and $\gamma$ arbitrary, $\kappa=0$. If in
(\icinl{ii}), (\icinl{iv}) 
and (\icinl{v}) $c_4=c_5=0$, $c_3\not=0$, then we may set $c_2=0$ and
$c_3=1$ without loss of 
generality. Thus we obtain the reduction
\[
 u(x,t) = t^{-2} w(z) , \qquad z=x- c_1 \log (t) ,
\]
where $w(z)$ satisf\/ies
\[
\ba{l}
\ds w \wzzzz+ \alpha \wz \wzzz + \beta \left(\wzz\right)^2
\vspace{3mm}\\
\ds \qquad + 2 \gamma \left[w \wzz + \left(\wz\right)^2\right] -c_1^2
\wzz -5c_1 \wz -6w=0.
\ea
\]

\redn{iciii} $\alpha$, $\beta$ and $\kappa$ arbitrary, $\gamma=0$. If in
(\icinl{iii}) $c_3\not=0$, then we may set $c_1=c_2=0$ and $c_3=1$
without loss of generality. Thus we obtain the reduction
\[
 u(x,t) =t^2 w(z) , \qquad z=x / t,
\]
where $w(z)$ satisf\/ies
\[
w \wzzzz+ \alpha \wz \wzzz + \beta \left(\wzz\right)^2
+ \kappa \wzz -z^2 \wzz +2z \wz -2 w =0 .
\]

\redn{iciv} $\alpha$ and $\beta$ arbitrary, $\kappa=\gamma=0$. If in
(\icinl{iv}) and (\icinl{v}) $c_3=c_5=0$ and $c_4\not=0$, 
then we may set $c_1=0$ and $c_2=1$ without loss of generality. Thus we
obtain the reduction 
\[
 u(x,t) = w(z) \exp (4 c_4 t), \qquad z= x \exp (-c_4 t) ,
\]
where $w(z)$ satisf\/ies
\[
w \wzzzz+ \alpha \wz \wzzz + \beta \left(\wzz\right)^2
 -c_4^2 z^2 \wzz + 7 c_4^2 z \wz -16 c_4^2 w = 0 .
\]

\redn{icv} $\alpha$ and $\beta$ arbitrary, $\kappa=\gamma=0$. If in
(\icinl{iv}) and (\icinl{v}) $c_5=0$ and $c_3 c_4 \not=0$, 
then we may set $c_1=c_2=0$ and $c_3 =1$ without loss of generality. Thus we
obtain the reduction 
\[
u(x,t) = w(z) t^{4 c_4 -2}, \qquad z=x t^{-c_4},
\]
where $w(z)$ satisf\/ies
\[
\ba{l}
\ds w \wzzzz+ \alpha \wz \wzzz + \beta \left(\wzz\right)^2 
\vspace{3mm}\\
\ds \qquad  - c_4^2 z^2 \wzz 
+\left(7c_4^2 -5 c_4 \right)z \wz - \left( 16 c_4^2-20 c_4 +6\right) w
=0 . 
\ea
\]

\redn{icvi} $\alpha=-\tfr52$, $\beta=\tfr{45}{28}$, $\gamma=\kappa=0$.
If in (\icinl{v}) $c_3=0$ and $c_5\not=0$, then we may set $c_1=-m c$, 
$c_2=1$, $c_4=0$ and $c_5 = m/c$,
without loss of generality. Thus we obtain the reduction
\[
 u(x,t) = {w(z) \exp (-8 m t) \over
 \left[z - \exp (-2 m t) \right]^8 }, \qquad
 z=\left( {x-c \over x+c} \right) \exp (-2 m t),
\]
where $w(z)$ satisf\/ies
\[
\ba{l}
\ds 28 w \wzzzz -70 \wz \wzzz +45 \left(\wzz\right)^2 
\vspace{3mm}\\
\ds \qquad -c^4 m^2 \left(1792 z^2 \wzz -12544 z \wz +28672 w\right)
=0. 
\ea
\]

\redn{icvii} $\alpha=-\tfr52$, $\beta=\tfr{45}{28}$, $\gamma=\kappa=0$.
If in (\icinl{v}) $c_3 c_5\not=0$, then we may set $c_1=-m c$, $c_2=
c_4=0$, $c_3=1$ and $c_5 =m/c$, without loss of 
generality. Thus we obtain the reduction
\[
u(x,t) = {w(z) t^{-2(1+4 m)} \over (z- t^{-2m})^8 }, \qquad 
 z=\left( {x-c \over x+c} \right) t^{-2 m},
\]
where $w(z)$ satisf\/ies
\[
\ba{l}
\ds 28 w \wzzzz -70 \wz \wzzz +45 \left(\wzz\right)^2 
 - 1792 c^4 m^2 z^2 \wzz
\vspace{3mm}\\ 
\ds \qquad + \left(12544 m^2-4480 m\right)c^4 z \wz
 -\left( 28672 m^2 -17920 m +2688 \right) c^4w  =0. 
\ea
\]

\subsection{Case (ii) $\mu\not=0$.}
In this case we generate 18 \deqs,
\[
\ba{l}
\tau_{u}=0, \qquad \tau_{x}=0, \qquad
 \xi_{u}=0, \qquad \xi_{t}=0, \qquad \phi_{uu}=0,
\vspace{1mm}\\
 \phi_{xtu}=0, \qquad \alpha (\phi_{u} u - \phi)=0, \qquad 2 \phi_{tu} -
\tau_{tt}=0, \qquad \beta (\phi_{u} u - \phi)=0,
\vspace{1mm}\\
2 \phi_{xu} - \xi_{xx}=0, \qquad 4 \phi_{xu} u - 6 \xi_{xx} u + \alpha
\phi_{x}=0,
\vspace{1mm}\\ 
2 \tau_{t} u - 2 \xi_{x} u + \phi=0,
\vspace{1mm}\\
4 \beta \phi_{xu} + 3 \alpha \phi_{xu} - 2 \beta \xi_{xx} - 3 \alpha
\xi_{xx}=0,
\vspace{1mm}\\ 
\phi_{xxu} u + 2 \tau_{t} u - 4 \xi_{x} u + \phi=0,
\vspace{1mm}\\
3 \alpha \phi_{xxu} u + 2 \gamma \phi_{u} u + 4 \xi_{x} \gamma u - \alpha
\xi_{xxx} u - 2 \gamma \phi=0,
\vspace{1mm}\\
\phi_{tt}- \phi_{xxxx} u - 2 \gamma \phi_{xx} u - \kappa \phi_{xx} -
\phi_{xxtt}=0,
\vspace{1mm}\\ 
6 \phi_{xxu} u^{2} + 4 \xi_{x} \gamma u^{2} - 4 \xi_{xxx} u^{2} + 2 \beta
\phi_{xx} u + \phi_{ttu} u  + 2 \xi_{x} \kappa u - \kappa \phi=0,
\vspace{1mm}\\
4 \phi_{xxxu} u + 4 \gamma \phi_{xu} u - 2 \xi_{xx} \gamma u - \xi_{xxxx} u
 + \alpha \phi_{xxx} + 4 \gamma \phi_{x} 
+ 2 \kappa \phi_{xu} - \xi_{xx} \kappa=0,
\ea
\]
and then use {\tt reduceall} in {\tt diffgrob2} to simplify them 
to the following system,
\[
\xi_{u}=0, \quad \xi_{t}=0, \quad
\xi_{x}=0, \quad \tau_{u} =0, \quad \kappa \tau_{t}=0, \quad
\tau_{x} =0, \quad \kappa \phi =0.
\]
Here $\kappa=0$ is the only special case, yielding the slightly dif\/ferent
system 
\[
\xi_{u}=0, \quad \xi_{t}=0, \quad
\xi_{x}=0, \quad \tau_{u} =0, \quad \tau_{tt}=0, \quad
\tau_{x} =0, \quad \phi+2 \tau_t u  =0. 
\]
Thus we have two dif\/ferent sets of \inls, and in both cases $\alpha, \beta$
and $\gamma$ remain arbitrary.\def\iicinl{2.2}

\begin{center}
{\bf Table \iicinl}

\medskip
\begin{tabular}{|ccccr|}
\hline & & & &\\[-5pt]
Parameters & $\xi$ & $\tau$ & $\phi$ & \\[5pt] \hline  & & & &\\[-5pt]
& $c_1$ & $c_2$ & $0$ & \eqtab{\iicinl}{i} \\[5pt] 
$\kappa=0$ & $c_1$ & $c_3 t +c_2$ & $-2 c_3 u$ &\eqtab{\iicinl}{ii} \\[5pt]
\hline 
\end{tabular}
\end{center}

From these we have the following two canonical reductions:

\redn{iici} $\alpha$, $\beta$, $\gamma$ and $\kappa$ arbitrary. If
in (\iicinl{i}) and (\iicinl{ii}) $c_3=0$, then we may set $c_2=1$ 
without loss of generality. Thus we obtain the 
following reduction
\[
u(x,t) = w(z), \qquad z= x-c_1 t,
\]
where $w(z)$ satisf\/ies
\[
\ba{l}
\ds \left( \kappa -c_1^2\right) \wzz + 2 \gamma \left[w \wzz +
\left(\wz\right)^2\right]
\vspace{3mm}\\
\ds \qquad + w \wzzzz+ \alpha \wz \wzzz + \beta
\left(\wzz\right)^2 +c_1^2 \wzzzz=0 .
\ea
\]

\redn{iicii} $\alpha$, $\beta$, $\gamma$ arbitrary, $\kappa=0$.
If in (\iicinl{ii}) $c_3\not=0$, then we may set $c_2=0$, $c_3=1$
without loss of generality. Thus we obtain the following reduction
\[
u(x,t) = t^{-2} w(z), \qquad  z= x- c_1 \log (t),
\]
where $w(z)$ satisf\/ies
\[
\ba{l}
\ds w \wzzzz+ \alpha \wz \wzzz + \beta \left(\wzz\right)^2 + c_1^2
\wzzzz 
\vspace{3mm}\\ 
\ds \qquad 
+5 c_1 \wzzz+ 2 \gamma \left[w \wzz + \left(\wz\right)^2\right] 
+ \left(6 -c_1^2\right) \wzz
- 5 c_1 \wz - 6 w =0 .
\ea
\]

\subsection{Travelling wave reductions}
As was seen in \S~1, special cases of (\ref{fulleqn}) admit interesting
travelling wave solutions, namely compactons. In this subsection we look for
such solitary waves and others, in the framework of (\ref{fulleqn}).
Starting with compacton-type solutions, we seek solutions of the form
\beq 
u(x,t) = a_2 \cos^n \{ a_3 (x-a_1 t) \} + a_4, \label{cosn} 
\eeq
where $a_1,a_2,a_3,a_4$ are constants to be determined. We include the
(possibly non-zero) constant $a_4$ since $u$ is open to translation. The
specif\/ic form of the translation will put conditions on $a_4$, which may or
may not put further conditions on the other parameters in (\ref{cosn}) and
those in (\ref{fulleqn}) -- see below. If 
$n=1$ we have the solutions, where the absence of a parameter implies it is
arbitrary, 
\[
\ba{rl}
(\mbox{i}) & \ds \alpha=0,\quad \beta=-1,\quad \gamma=0,
\quad a_4 = {\kappa - a_1^2 - a_1^2 a_2^2 \mu \over a_2^2} .
\vspace{3mm}\\
(\mbox{ii})  & \ds \alpha=1, \quad\beta=0,\quad \gamma>0 ,\quad
a_1^2(1+2 \gamma \mu) - \kappa =0 ,\quad a_3^2 = 2\gamma.
\vspace{3mm}\\
(\mbox{iii}) & \ds \beta=\alpha-1\not=0,\quad {\gamma \over \alpha}>0,
\quad a_3^2 = {2\gamma \over \alpha},\quad a_4 = {\alpha(\kappa - a_1^2) -2
a_1^2 \gamma\mu \over 2 \gamma (1-\alpha) }.
\ea
\]
These become $n=2$ solutions via the trigonometric identity $\cos 2 \theta =
2 \cos^2 \theta -1$. By earlier reasoning the associated compactons are weak
solutions of (\ref{fulleqn}). When considering more general $n$ we restrict
$n$ to be either 3 or $\ge 4$ else the fourth derivatives of $u(x,t)$ that
we require in (\ref{fulleqn}) would have singularities at the edges of the
humps; we f\/ind 
\[
\alpha={2\over n}, \quad\beta = {2-n \over n},\quad\gamma>0 ,
\quad  a_1^2(1+2 \gamma \mu) - \kappa =0 ,\quad a_3^2 = {\gamma \over n} ,\quad
a_4 = -\mu a_1^2 .
\]
When $n=3$ or $n=4$ our compacton would be a weak solution since not all the
derivatives of $u(x,t)$ in (\ref{fulleqn}) in these instances are continuous
at the edges. For $n>4$ the solutions are strong.

For more usual solitary waves we seek solutions of the form
\[
u(x,t) = a_2 \mbox{sech}^n \{ a_3 (x-a_1 t) \} + a_4, \label{sechn} 
\]
where $a_1,a_2,a_3,a_4$ are constants to be determined. If $n=2$ then
$\alpha=-1$, $\beta=-2$ and we have solutions
\[
\ba{rl}
(\mbox{i})  &\ds \gamma<0 ,\quad a_1^2(1+2 \gamma \mu) - \kappa =0 ,\quad
a_3^2 = -\tfr12\gamma ,\quad a_4 = -\tfr13 ( 2a_2 +3a_1^2 \mu),
\vspace{3mm}\\
(\mbox{ii}) & \ds \gamma^2\not=4a_3^4, \quad a_2 = { 3 a_3^2
( \kappa - a_1^2 -2a_1^2\gamma\mu)\over (2a_3^2-\gamma)(2a_3^2+\gamma)},
\quad a_4 = -{ \kappa - a_1^2+4a_1^2 a_3^2 \mu \over 2(2a_3^2+\gamma)} ,
\ea
\]
and for general $n$, including $n=2$ ($\gamma>0$)
\[
\alpha=-{2\over n}, \qquad\beta = -{n+2\over n},\qquad
 a_1^2(1+2 \gamma \mu) - \kappa =0 ,\qquad
a_3^2 = {\gamma \over n} ,\qquad
a_4 = -\mu a_1^2 .
\]
Now consider the general travelling wave reduction, $u(x,t)=w(z)$,
$z=x-ct$, where $w(z)$ satisf\/ies
\[
\ba{l}
\ds \left(\kappa- c^2\right) \wzz +2\gamma
\left[w\wzz+\left(\wz\right)^2\right]
\vspace{3mm}\\
\ds \qquad +\mu c^2 \wzzzz+ w\wzzzz +\alpha \wz\wzzz +
\beta \left(\wzz\right)^2 =0.
\ea
\]
In the special case $\beta=\alpha-1$, we can integrate
this twice with respect to $z$ to give
\beq 
\left(\kappa- c^2\right) w+\gamma w^2 +\mu c^2 \wzz +\frac 12 (\alpha-2)
\left(\wz\right)^2 +A z +B =0 \label{intwice}
\eeq
with $A$ and $B$ the constants of integration.
If we assume $A=0$, then we make the transformation $W(z)=w(z)+\mu c^2$,
multiply (\ref{intwice}) by $\displaystyle W^{\alpha-3} \Wz $ and integrate
\wrt\ $z$ to yield
\[
{\gamma \over \alpha} W^{\alpha} + {A_1 \over \alpha-1}
W^{\alpha-1} + {A_2 \over \alpha-2} W^{\alpha-2} +
 W^{\alpha-2} \left(\Wz\right)^2 + C=0 \label{Weqn}
\]
for $\alpha\not=0,1,2$, where $A_1= \kappa - c^2 -2\gamma \mu c^2$
and $A_2 = B - \mu c^2 \left(\kappa - c^2 -\gamma \mu c^2\right)$.
In the special cases $\alpha=0,1,2$ we obtain respectively
\be
\gamma\log W - {A_1 \over W} - {2A_2 \over W^2} + {1 \over
W^2}\left(\Wz\right)^2 +C=0, \label{alpWeqni} 
\ee
\be
\gamma W + A_1 \log W - {A_2 \over W} + {1 \over W}\left(\Wz\right)^2 +C=0,
\label{alpWeqnii} 
\ee
\be
\frac 12\gamma W^2 +A_1 W + A_2 \log W + \left(\Wz\right)^2 + C=0 ,
\label{alpWeqniii}
\ee
where $C$ is a constant of integration. For $\alpha\in{\mathbb Z}$,
an integer, with $\alpha\ge 3$, (\ref{Weqn}) may be written as
\be
\ba{l}
\ds
\ds \left(w+\mu c^2\right)^{\alpha-2} \left(\wz\right)^2 
\vspace{3mm}\\
\ds \qquad +{\gamma w^2\over\alpha} 
 \left[ \left(w+\mu c^2\right)^{\alpha-2} + { \alpha 
\left(\kappa - c^2\right) - 
2\mu \gamma c^2\over\gamma} \sum\limits_{n=0}^{\alpha-3}
{(\mu c^2)^{\alpha-3-n}\over(n+2)}
{\alpha-3 \choose n} {w^n}\right] 
\vspace{3mm}\\ 
\ds \qquad +{B \over \alpha-2}
\left(w+\mu c^2\right)^{\alpha-2} + 
C - {(\mu c^2)^{\alpha-1} \left[ \alpha \left(\kappa - c^2\right) - 2\mu
\gamma c^2\right] \over \alpha (\alpha-1) (\alpha-2) }=0,
\ea
\label{ppeakeq}
\ee
where $C$ is a constant of integration.
If we require that $w$ and its derivatives tend to zero as $z\rightarrow
\pm\infty$, then $B=D=0$. If $\alpha=3$ this equation induces so-called
peakons (cf.\ 
\cite{refCH}) as $\alpha \left(\kappa - c^2\right) - 2\mu \gamma
c^2\rightarrow 0$ (see \cite{refCHH,refGP,refK,refRosa}). Similarly if
$\alpha=4$ this equation is of the form found in \cite{refGP} which induces
the `wave of greatest height' found in \cite{refFW}. Both solutions, in
their limit, have a discontinuity in their f\/irst derivative at its peak.
Note that if $\alpha \left(\kappa - c^2\right) - 2\mu \gamma c^2=0$, 
equation (\ref{ppeakeq}) becomes
\beq 
\left(w+\mu c^2\right)^{\alpha-2} \left[ \left(\wz\right)^2 + {\gamma
\over\alpha} w^2 \right] =0. 
\label{peakeqn}
\eeq
Since $\alpha>0$ then we require $\gamma<0$ to give a peakon of the form
\beq 
u(x,t) = {\alpha \left(c^2 -\kappa\right)\over 2\gamma} \exp
\left\{ - \left( - \gamma\over\alpha \right)^{1/2} | x-ct| \right\}.
\label{peak} 
\eeq
The height of the wave, because of the form of (\ref{fulleqn}), is dependent
upon the square of the speed, whereas the peakons in \cite{refCH} and
\cite{refFW} are proportional to the wave speed.

\def\sen{3}
\redno=0

\section{Nonclassical symmetries ($\tau\not=0$)}
In the nonclassical method one requires only the subset of $S_{\Delta}$
given by 
\beq 
S_{\Delta,\psi} = \{ u(x,t) : \Delta (u) =0, \psi (u) =0 \}, \label{sdelpsi}
\eeq 
where $S_{\Delta}$ is def\/ined in (\ref{sdel}) and $\psi=0$ is the \isc\
(\ref{invsc}), to be invariant under the transformation (\ref{trans}). The
usual method of applying the nonclassical method (e.g. as described in
\cite{refLW}), involves applying the prolongation ${\rm pr}^{(4)} {\bf v}$
to the system composed of (\ref{fulleqn})\ and the \isc\ (\ref{invsc}) and
requiring that the resulting expressions vanish for $u\in S_{\Delta,\psi}$,
i.e. 
\beq
{\rm pr}^{(4)} {\bf v} (\Delta) |_{\{\Delta=0,\psi=0\}}=0,
\qquad {\rm pr}^{(1)} {\bf v} (\psi) \vert_{\{\Delta=0,\psi=0\}}=0.
\label{twoprol}
\eeq
It is well known that the latter vanishes identically when $\psi=0$ without
imposing any conditions upon $\xi$, $\tau$ and $\phi$. To apply the method
in practice we advocate the algorithm described in \cite{refCMiii} for
calculating the determining equations, which avoids dif\/f\/iculties arising
from using dif\/ferential consequences of the \isc\ (\ref{invsc}). 

In the canonical case when $\tau\not=0$ we set $\tau=1$ without loss of
generality. We proceed by eliminating $u_{tt}$ and $u_{xxtt}$ in
(\ref{fulleqn})\ using the \isc\ (\ref{invsc}) which yields
\be
\ba{l}
 \xi \xi_{x}u_x + 2 u_x^{2} \xi \xi_{u}- 2 \phi_{u} \xi u_x
+ \xi^{2}u_{xx}- \phi_{x} \xi - \xi_{t}u_x 
+ \phi \phi_{u}- \phi \xi_{u} u_x+ \phi_{t}
 - \kappa u_{xx}
\vspace{1mm} \\  \quad
   - 2 \gamma ( u u_{xx} + u_x^{2})
- u u_{xxxx} - \alpha u_x u_{xxx} - \beta u_{xx}^{2}
+ \mu \Bigl[ 2 \phi_{xx} \xi_{x}- 2  \phi_{xu} \phi_{x}
- 4   \xi_{x}^{2} u_{xx}
\vspace{1mm}\\ \quad
      -  \phi_{tu} u_{xx} -  \phi_{u} \phi_{xx}
      -  \phi_{xxt} +  \xi_{xxt} u_x -  \phi_{tuu} u_x^{2}
      +  \phi_{x} \xi_{xx} +  \xi_{tuu} u_x^3 
      +  \xi_{t} u_{xxx} -  \xi^{2} u_{xxxx}
\vspace{1mm}\\ \quad
       +  \phi_{xxx} \xi
      -  \phi_{u}^{2} u_{xx} -  \phi \phi_{xxu}
      -  \xi \xi_{xxx} u_x +  \phi \xi_{uuu} u_x^3
      +  \phi \xi_{xxu} u_x +  \phi \xi_{u} u_{xxx}
      -  \phi \phi_{uu} u_{xx}
\vspace{1mm}\\ \quad
       -  \phi \phi_{uuu} u_x^{2}
      + 2  \xi_{xt} u_{xx} + 2  \xi_{xtu} u_x^{2} 
      - 2  \phi_{xtu} u_x - 3  \xi_{x} \xi_{xx} u_x
      - 4  \xi_{u} \xi_{xx} u_x^{2}
      - 4  \xi \xi_{xx} u_{xx}
\vspace{1mm}\\ \quad
      + 2  \phi_{u} \xi_{xx} u_x- 5  \xi_{uu} \xi_{x} u_x^{3}
      - 8  \xi_{xu} \xi_{x} u_x^{2}
      - 15 \xi_{u} \xi_{x} u_x u_{xx} 
      - 5  \xi \xi_{x} u_{xxx}
      + 4  \phi_{u} \xi_{x} u_{xx}
\vspace{1mm}\\ \quad
      + 4  \phi_{uu} \xi_{x} u_x^{2}
      + 6  \phi_{xu} \xi_{x} u_x
      - 2  \xi \xi_{uuu} u_x^{4}- 5  \xi \xi_{xuu} u_x^{3}
      + 2  \phi \xi_{xuu} u_x^{2}- 6  \xi_{u} \xi_{uu} u_x^{4}
\vspace{1mm}\\ \quad
      - 12  \xi \xi_{uu} u_x^{2} u_{xx}
      + 3  \phi \xi_{uu}u_x u_{xx} 
      + 4  \phi_{u} \xi_{uu} u_x^{3}+ 3  \phi_{x} \xi_{uu} u_x^{2}
      - 4  \xi \xi_{xxu} u_x^{2}- 10  \xi_{u} \xi_{xu} u_x^{3}
\vspace{1mm}\\ \quad
      - 15  \xi \xi_{xu} u_x u_{xx}
      + 2  \phi \xi_{xu}u_{xx} + 6  \phi_{u} \xi_{xu} u_x^{2} 
     + 4  \phi_{x} \xi_{xu} u_x
      - 12  \xi_{u}^{2} u_x^{2} u_{xx}
      - 8  \xi \xi_{u} u_x u_{xxx}
\vspace{1mm}\\ \quad
      - 6  \xi \xi_{u} u_{xx}^{2}
      + 9  \phi_{u} \xi_{u} u_x u_{xx}
      + 3  \phi_{x} \xi_{u} u_{xx}+ 5  \phi_{uu} \xi_{u} u_x^{3} 
      + 8  \phi_{xu} \xi_{u} u_x^{2}+ 3  \phi_{xx} \xi_{u} u_x
\vspace{1mm}\\ \quad
      + 3  \xi_{tu} u_x u_{xx}
      + 2  \phi_{u} \xi u_{xxx}
      + 6  \phi_{uu} \xi u_x u_{xx}
      + 5  \phi_{xu} \xi u_{xx}+ 2  \phi_{uuu} \xi u_x^{3} 
      + 5  \phi_{xuu} \xi u_x^{2}
\vspace{1mm}\\ \quad
      + 4  \phi_{xxu} \xi u_x
      - 3  \phi_{u} \phi_{uu} u_x^{2} - 2  \phi_{uu} \phi_{x} u_x
      - 2  \phi \phi_{xuu} u_x
      - 4  \phi_{u} \phi_{xu} u_x\Bigr] =0. \label{delstar}
\ea \hspace{-14pt}
\ee
We note that this equation now involves the inf\/initesimals $\xi$ and $\phi$
that are to be determined. Then we apply the classical Lie algorithm
to (\ref{delstar}) using the fourth prolongation $ {\rm pr}^{(4)}
{\bf v}$ and eliminating $u_{xxxx}$ using (\ref{delstar}). It should
be noted that the coef\/f\/icient of $u_{xxxx}$ is ($\xi^2+\mu u$).
Therefore, if this is zero the removal of $u_{xxxx}$ using
(\ref{delstar}) is invalid and so the next highest derivative term,
$u_{xxx}$, should be used instead. We note again that this has a
coef\/f\/icient that may be zero so that in the case $\mu\not=0$ and
$\xi^2+\mu u=0$ one needs to calculate the \deqs\ for the subcases
non-zero separately. Continuing in this fashion, there is a cascade of 
subcases to be considered.  In the remainder of this section, we
consider these subcases in turn.  First, however, we discuss
the case given by~$\mu=0$.

\subsection{Case (i) $\mu=0$.}
In this case we generate the following 12 \deqs.
\[
\ba{l}
\xi_{u}=0,
\vspace{1mm}\\
 \phi_{uuuu} u + \alpha \phi_{uuu}=0,
\vspace{1mm} \\
4 \phi_{xuuu} u + 3 \alpha \phi_{xuu} =0,
\vspace{1mm}\\
6 \phi_{uuu} u + 2 \beta \phi_{uu} + 3 \alpha \phi_{uu}=0,
\vspace{1mm}\\
4 \phi_{uu} u^{2} + \alpha \phi_{u} u - \alpha \phi=0,
\vspace{1mm}\\
4 \phi_{xu} u - 6 \xi_{xx} u + \alpha \phi_{x}=0,
\vspace{1mm}\\
3 \phi_{uu} u^{2} + \beta \phi_{u} u - \beta \phi=0,
\vspace{1mm}\\
12 \phi_{xuu} u + 4 \beta \phi_{xu} + 3 \alpha \phi_{xu} - 2 \beta \xi_{xx}
        - 3 \alpha \xi_{xx}=0,
\vspace{1mm} \\
6 \phi_{xxuu} u^{2} + 2 \gamma \phi_{uu} u^{2} + \kappa \phi_{uu} u
          - \xi^{2} \phi_{uu} u
\vspace{1mm}\\
\qquad  + 3 \alpha \phi_{xxu} u + 2 \gamma \phi_{u} u
          + 4 \xi_{x} \gamma u - \alpha \xi_{xxx} u - 2 \gamma \phi=0,
\vspace{1mm}\\
6 \phi_{xxu} u^{2} + 4 \xi_{x} \gamma u^{2} - 4 \xi_{xxx} u^{2}
          + 2 \beta \phi_{xx} u + 2 \xi_{x} \kappa u - 4 \xi^{2} \xi_{x} u
          - 2 \xi \xi_{t} u - \kappa \phi + \xi^{2} \phi =0,
\vspace{1mm}\\
\phi_{tt} u - \phi_{xxxx} u^{2} -
2 \gamma \phi_{xx} u^{2} - \kappa \phi_{xx} u
- 4 \xi \xi_{x} \phi_{x} u - 2 \xi_{t} \phi_{x} u
\vspace{1mm}\\
\qquad  + \phi^{2} \phi_{uu} u + 4 \xi_{x} \phi \phi_{u} u
        + 2 \phi \phi_{tu} u + 4 \xi_{x} \phi_{t} u + \xi \phi \phi_{x}
        - \phi^{2} \phi_{u} - \phi \phi_{t}=0,
\vspace{1mm}\\
4 \phi_{xxxu} u^{2} + 4 \gamma \phi_{xu} u^{2} - 2 \xi_{xx} \gamma u^{2}
          - \xi_{xxxx} u^{2}
\vspace{1mm}\\
\qquad + \alpha \phi_{xxx} u + 4 \gamma \phi_{x} u
          + 2 \xi \phi \phi_{uu} u + 2 \kappa \phi_{xu} u
          + 8 \xi \xi_{x} \phi_{u} u + 2 \xi_{t} \phi_{u} u
\vspace{1mm}\\
\qquad   + 2 \xi \phi_{tu} u - \xi_{xx} \kappa u - 4 \xi \xi_{x}^{2} u
          + 2 \xi_{t} \xi_{x} u + \xi_{tt} u - 2 \xi \phi \phi_{u}
          + \xi \xi_{x} \phi - \xi_{t} \phi=0.
\ea
\]
As guaranteed by the nonclassical method, we get all
the classical reductions, but we also have some
\inls\ that lead to nonclassical reductions, namely

\def\ininl{3.1}
\begin{center}
{\bf Table \ininl}

\medskip
\begin{tabular}{|cccr|}
\hline & & &\\[-5pt]
Parameters & $\xi$ & $\phi$ &  \\[5pt] \hline & & &\\[-5pt]
$\kappa=0$ & $0$  & $g(t) u$
\ where  $\displaystyle\gtt + g \gt -g^3 =0$ & \eqtab{\ininl}{i} \\[5pt]
$\alpha=\beta=\gamma=0$ & $\pm \sqrt{\kappa}$ & 
$c_3 y^3 +c_2 y^2 + c_1 y + c_0 \quad (y=x \pm \sqrt{\kappa}\, t)$ 
&\eqtab{\ininl}{ii}  \\[5pt]
$\alpha=\beta=\gamma=\kappa=0$ & $0$ &
$\ba{l} \ds -  \frac{u}{g(t)}\gt + g(t)\left( c_4 x^4 +c_3 x^3 +c_2 x^2 +
c_1 x + c_0\right)\\[10pt]
\mbox{where}\\[5pt]
 \displaystyle g^2 \gttt -4 g \gt \gtt
+ 2 \left(\gt\right)^3 + 24 c_4 g^4 =0 \ea $  
& \eqtab{\ininl}{iii} \\[10pt]   \hline
\end{tabular}
\end{center}

From these we obtain three canonical reductions.

\redn{ini} $\alpha$, $\beta$, $\gamma$ arbitrary, $\kappa=0$. In
(\ininl{i}) we solve the equation for $g(t)$ by writing
$g(t)=[\ln(\psi(t))]_t$ then $\psi(t)$ satisf\/ies
\beq
\left(\frac{\d\psi}{\d t}\right)^2 = 4 c_1 \psi^3 + c_2, \label{psieqn}
\eeq
where $c_1$ and $c_2$ are arbitrary constants;
$c_1=c_2=0$ is not allowed 
since $g(t)\not\equiv0$.
Hence we obtain the following reduction
\[
 u(x,t) = w(x) \psi (t),
 \]
where $w(x)$ satisf\/ies
\[
w \wxxxx+ \alpha \wx \wxxx + \beta \left(\wxx\right)^2
+ 2 \gamma \left[w \wxx + \left(\wx\right)^2\right] - 6 c_1 w =0.
\]
There are three cases to consider in the solution of (\ref{psieqn}). 

\begin{itemize}
\item[(i)] If $c_1=0$, we may assume that $\psi(t) = t$
without loss of generality.

\item[(ii)] If $c_2=0$, then $\psi= \left[ c_2 (t+c_3)^2\right]^{-1}$
and we may set $c_2=1$, $c_3=0$ without loss of generality.

\item[(iii)] If $c_1 c_2 \not=0$ we may set $c_1=1$,
$c_2=-g_3$ without loss of generality so that $\psi(t)$
is any solution of the  Weierstrass elliptic function
equation
\beq
\left[\frac{\d\wp}{\d t}(t;0,g_3)\right]^2 =
4 \wp^3 (t;0,g_3) - g_3 .\label{ellip}
\eeq

\end{itemize}

\redn{inii} $\kappa$ arbitrary, $\alpha=\beta=\gamma=0$. From
(\ininl{ii}) we get the following reduction
\[
 u(x,t) = w(z) \pm {c_3 \over 8 \sqrt{\kappa}} y^4 \pm
 {c_2 \over 6 \sqrt{\kappa}} y^3\pm
 {c_1 \over 4 \sqrt{\kappa}} y^2 + c_0 t, \! \! \!
 \qquad y=x \pm \sqrt{\kappa}\, t,
\quad  z= x \mp \sqrt{\kappa}\, t,
\]
where $w(z)$ satisf\/ies
\[
\sqrt{\kappa} \,\wxxxx \pm 3 c_3 =0.
\]
This gives us the exact solution
\[
u(x,t) = \mp {c_3 \over 8 \sqrt{\kappa}} z^4 + c_4 z^3+
c_5 z^2 + c_6 z + c_7 \pm {c_3 \over 8 \sqrt{\kappa}} y^4 \pm
{c_2 \over 6 \sqrt{\kappa}} y^3\pm
{c_1 \over 4 \sqrt{\kappa}} y^2 + c_0 t.
\]

\redn{iniii} $\alpha=\beta=\gamma=\kappa=0$.
In (\ininl{iii}) we integrate our
equation for $g(t)$ up to an expression with quadratures
\beq
g \gtt - 2 \left(\gt\right)^2 + 24 c_4 g \int^t g^2 (s) \, \d s
+24 c_5 g = 0 .\label{geqn}
\eeq
We get the following reduction
\[
u(x,t) =  \frac{1}{g(t)} \left[w(x) +
\left(c_4 x^4 +c_3 x^3 +c_2 x^2 + c_1 x + c_0\right)
\int^t g^2 (s) \, \d s\right],
\]
where $w(x)$ satisf\/ies
\[
\wxxxx - 24 c_5 = 0.
\]
This is easily solved to give the solution
\[
\ba{l}
\ds u(x,t) =  \frac{1}{g(t)} \Biggl[ c_5 x^4 + c_6 x^3 +c_7 x^2 + c_8 x + c_9
\vspace{3mm}\\
\ds \qquad + \left(c_4 x^4 +c_3 x^3 +c_2 x^2 + c_1 x + c_0\right)
\int^t g^2 (s) \, \d s\Biggr],
\ea
\]
where $g(t)$ satisf\/ies (\ref{geqn}).

\subsection{Case (ii) $\mu\not=0$.} 
As discussed earlier in this section, we must consider,
in addition to the general case of the determining equations,
each of the singular cases of the determining equations.

\subsubsection{\mbox{$\xi^2+u\not=0$.}}
In this the generic case we generate 12 \deqs\ -- see appendix A
for details. As expected we have all the classical reductions,
however we also have the following \inls\ that
lead to genuine nonclassical reductions (i.e.\ not a classical reduction).

\def\iininl{3.2}
\begin{center}
{\bf Table \iininl}

\medskip
\begin{tabular}{|cccr|}
\hline & & & \\[-5pt]
Parameters & $\xi$ & $\phi$ &  \\[5pt]
\hline & & & \\[-5pt]
$\kappa=0$ & $0$  & $g(t) u$ \ where
$\displaystyle\gtt + g \gt -g^3 =0$ &  \eqtab{\iininl}{i}\\[5pt]
$1+2 \gamma= 0$ & $c_1 t+c_2$ &
$-2 c_1 (c_1 t+c_2)$ & \eqtab{\iininl}{ii}\\[5pt]
$\kappa=1+2 \gamma= 0$ & $c_2 ( t+c_1)^2$ &
$u (t+c_1)^{-1} - 3 c_2^2 (t+c_1)^3$ & \eqtab{\iininl}{iii}\\[5pt]
$\alpha=\beta=\gamma=0$ & $\pm \sqrt{\kappa}$ &
$\!\!c_3 y^3 +c_2 y^2 + c_1 y + c_0 \ (y=x \pm \sqrt{\kappa} t)\!\!$
& \eqtab{\iininl}{iv}\\[5pt]
$\alpha=-\tfr32, \beta=2, \gamma=0$ & $\pm \tfr12\sqrt{\kappa}\, (x+c_1)$
& $\pm 2 \sqrt{\kappa}\, u \pm \tfr14\kappa^{3/2} 
 (x+c_1)^2$ & \eqtab{\iininl}{v} \\[5pt]
$\alpha=\beta=\gamma=\kappa=0$ & $0$ & $c_3 x^3 +c_2 x^2 + c_1 x
+ c_0$ & \eqtab{\iininl}{vi}\\[5pt]
$\alpha=\beta=\gamma=\kappa=0$ & $0$ & $\!\!(u+c_3 x^3 +c_2 x^2 + c_1 x
+ c_0)( t+c_4)^{-1}\!\!$ & \eqtab{\iininl}{vii}\\[10pt] \hline
\end{tabular}
\end{center}

From these inf\/initesimals we obtain six reductions.

\redn{iini} $\alpha$, $\beta$,
$\gamma$ arbitrary, $\kappa=0$. In (\iininl{i}) we solve
the equation for $g(t)$ by writing
$g(t)=[\ln(\psi(t))]_t$ then $\psi(t)$ satisf\/ies
\beq
\left(\frac{\d\psi}{\d t}\right)^2 = 4 c_1 \psi^3 + c_2 \label{psieqnn}
\eeq
though $c_1=c_2=0$ is not allowed to preserve the fact that
$g(t)\not\equiv0$. We obtain the following reduction
\[
 u(x,t) = w(x) \psi (t),
\]
where $w(x)$ satisf\/ies
\[
w \wxxxx+ \alpha \wx \wxxx + \beta \left(\wxx\right)^2
+ 2 \gamma \left[w \wxx + \left(\wx\right)^2\right]
+ 6 c_1\left(\wxx- w\right) =0 .
\]
There are three cases to consider in the solution of (\ref{psieqnn}). 

\begin{itemize}
\item[(i)] If $c_1=0$, we may assume that
$\psi(t) = t$ without loss of generality.

\item[(ii)] If $c_2=0$, then $\psi= \left[ c_2 (t+c_3)^2\right]^{-1}$
and we may set $c_2=1$ and $c_3=0$ without loss of generality.

\item[(iii)] If $c_1 c_2 \not=0$ we may set $c_1=1$ and $c_2=-g_3$
without loss of generality so that $\psi(t)$ is any solution of
the Weierstrass elliptic function equation (\ref{ellip}).
\end{itemize}

Note that in the special case
\[
\wzz-w=0,
\]
we are able to lift the restrictions on
$\psi(t)$ so that it is arbitrary, if
$\beta+1+2 \gamma = \alpha+2 \gamma =0$. This
yields the exact solution
\[
u(x,t)  = \psi (t) \left(c_2 \e^x +c_3 \e^{-x}\right),
\]
where $\psi (t)$ is arbitrary, $\kappa=0$, $\alpha=-2\gamma$
and $\beta = -1-2 \gamma$.

\redn{iinii} $\alpha$, $\beta$ and $\kappa$ are arbitrary and
$\gamma=-\tfr12$. In
(\iininl{ii}) we assume $c_1\not=0$ otherwise we get a
classical reduction, and may set $c_2=0$ without loss
of generality. Thus we obtain the following
accelerating wave reduction
\[
u(x,t) = w(z) - c_1^2 t^2 , \qquad z=x- \tfr12 c_1 t^2,
\]
where $w(z)$ satisf\/ies
\[
\ba{l}
\ds w \wzzzz+ \alpha \wz \wzzz + \beta \left(\wzz\right)^2
\vspace{3mm}\\
\ds \qquad  -c_1 \wzzz - w \wzz +\kappa \wzz - \left(\wz\right)^2 +c_1 \wz
 +2c_1^2 =0.
\ea
\]

\redn{iiniii} $\alpha$ and $\beta$ are arbitrary,
$\gamma=-\tfr12$ and $\kappa=0$. From
(\iininl{iii}) the following holds for arbitrary $c_2$,
and we may set $c_1=0$ without loss
of generality. Thus we obtain the reduction
\[
u(x,t) = w(z) t - c_2^2 t^4, \qquad z=x- \tfr13 c_2  t^3,
\]
where $w(z)$ satisf\/ies
\[
w \wzzzz+ \alpha \wz \wzzz + \beta \left(\wzz\right)^2
 -4 c_2 \wzzz -w \wzz - \left(\wz\right)^2 +4 c_2 \wz +12c_2^2 =0 .
\]

\redn{iiniv} $\kappa$ is arbitrary and $\alpha=\beta=\gamma=0$.
From (\iininl{iv}) we get the following reduction
\[
u(x,t) = w(z) \pm {c_3 \over 8 \sqrt{\kappa}} y^4 \pm
 {c_2 \over 6 \sqrt{\kappa}} y^3\pm
 {c_1 \over 4 \sqrt{\kappa}} y^2 + c_0 t, \! \! \!\qquad y=x \pm \sqrt{\kappa}\, t,
 \quad z= x \mp \sqrt{\kappa}\, t,
\]
where $w(z)$ satisf\/ies
\[
\sqrt{\kappa}\, \wzzzz \pm 3 c_3 =0 .
\]
This gives us the exact solution
\[
u(x,t) = \mp {c_3 \over 8 \sqrt{\kappa}} z^4 + c_4 z^3+
 c_5 z^2 + c_6 z + c_7 \pm {c_3 \over 8 \sqrt{\kappa}} y^4 \pm
 {c_2 \over 6 \sqrt{\kappa}} y^3\pm
 {c_1 \over 4 \sqrt{\kappa}} y^2 + c_0 t.
\]

\redn{iinv} $\kappa$ is arbitrary, $\alpha=-\tfr32$, $\beta=2$ and
$\gamma=0$. In (\iininl{v}) we may set $c_1=0$ without loss of generality.
Thus we obtain the following reduction
\[
 u(x,t) = w(z) x^4 - \tfr14\kappa x^2, \qquad
 z= \log (x) \mp \tfr12 \sqrt{\kappa}\, t,
\]
where $w(z)$ satisf\/ies
\[
\ba{l}
\ds 4w \wzzzz -6\wz \wzzz+8 \left(\wzz\right)^2
 + 16 w \wzzz +58 \wz \wzz
 \vspace{3mm}\\
 \ds \qquad +116 w \wzz - \kappa \wzz
+236 \left(\wz\right)^2 +776 w\wz +672 w^2 =0 .
\ea
\]

\redn{iinvi} $\alpha=\beta=\gamma=\kappa=0$. From (\iininl{vi}) and from
(\iininl{vii}) ($c_4=0$ without loss of generality) we get the
following reductions
\[
u(x,t) = w(x) + \left(c_3 x^3 +c_2 x^2 +c_1 x+c_0\right) t
\]
and
\[
u(x,t) = w(x) t - \left(c_3 x^3 +c_2 x^2 +c_1 x+c_0\right)
\]
respectively. In both cases $w(x)$ satisf\/ies
\[
\wxxxx = 0 .
\]
These reductions have a common exact solution, namely
\[
u(x,t) = P_3 (x) t + Q_3 (x),
\]
where $P_3$ and $Q_3$ are any third order polynomials in $x$.

\subsubsection{\mbox{$\xi^2+u=0$, not both $\alpha=4$ and
$2 \xi \phi_u+ \xi_u \phi=0$.}} 

The \deqs\ quickly lead us to
require that both $\alpha=4$ and $2 \xi \phi_u+
\xi_u \phi=0$, which is a contradiction. 

\subsubsection{\mbox{$\xi^2+u=0$, $\alpha=4$,
$2 \xi \phi_u+\xi_u \phi=0$ and $\beta\not=3$.}} 

The \deqs\ give us that $\gamma=-\tfr12$, $\kappa=0$ and $\phi=0$.
The \isc\ is then
\[
\pm \i \sqrt u\, u_x +u_t =0
\]
which may be solved implicitly to yield the solution
\[
u(x,t) = w(z), \qquad  z=x \mp \i \sqrt u\, t .
\]
However, substituting into our original equation
gives $\displaystyle\wz=0$, i.e.\ $u(x,t)$ is a constant.

\subsubsection{$\xi^2+u=0$, $\phi=H(x,t) u^{-1/4}$,
$\alpha=4$ and $\beta=3$, not all of
$H=0$, $\kappa=0$, $1+2 \gamma=0$.}

 For the \deqs\ to be satisf\/ied,
each of $H=0$, $\kappa=0$ and $1+2\gamma=0$, which is in
contradiction to our assumption.

\subsubsection{\mbox{$\xi^2+u=0$, $\phi=0$, $\alpha=4$, $\beta=3$,
$\gamma=-\tfr12$ and $\kappa=0$.}}

Under these conditions equation (\ref{delstar}) which we apply the
classical method to is identically zero. Therefore any solution of the
\isc\ is also a solution of (\ref{fulleqn}). Hence we get the following
reduction

\redn{iinvii} $\alpha=4$, $\beta=3$, $\gamma=-\tfr12$ and $\kappa =0$.
The \isc\ is
\[
\pm \i \sqrt u \,u_x +u_t =0
\]
which may be solved implicitly to yield
\[
u(x,t) = w(z), \qquad  z=x \mp \i \sqrt u\, t,
\]
where $w(z)$ is arbitrary.

\def\sen{4}
\redno=0

\section{Nonclassical symmetries ($\tau=0$)}

In the canonical case of the nonclassical method when $\tau=0$
we set $\xi=1$ without
loss of generality. We proceed by eliminating $u_x$, $u_{xx}$, $u_{xxx}$,
$u_{xxxx}$ and $u_{xxtt}$ in (\ref{fulleqn}) using the
\isc\ (\ref{invsc}) which yields
\be
\ba{l}
u_{tt} - \kappa (\phi_{x} + \phi \phi_{u})
 - 2 \gamma (u \phi_{x} + u \phi \phi_{u} + \phi^{2})
- u (\phi_{xxx} + \phi_{u} \phi_{xx} + \phi_{u}^{2} \phi_{x}
      + \phi \phi_{u}^{3}
\vspace{1mm}\\
\quad  + 4 \phi_{u} \phi^{2} \phi_{uu}
      + 5 \phi_{u} \phi \phi_{xu} + 3 \phi \phi_{uu} \phi_{x}
+ \phi^{3} \phi_{uuu}
+ 3 \phi^{2} \phi_{xuu} + 3 \phi \phi_{xxu}
      + 3 \phi_{xu} \phi_{x})
\vspace{1mm}\\
\quad
 - \mu (\phi \phi_{uu} u_{tt} + \phi \phi_{uuu} u_t^{2}
      + 2 \phi \phi_{tuu} u_t + \phi \phi_{ttu}
      + \phi_{u}^{2} u_{tt}
      + 3 \phi_{u} \phi_{uu} u_t^{2}
      + 4 \phi_{u} \phi_{tu} u_t
\vspace{1mm}\\
\quad + \phi_{u} \phi_{tt}
      + \phi_{xu} u_{tt} + \phi_{xuu} u_t^{2}
      + 2 \phi_{xtu} u_t + 2 \phi_{t} \phi_{uu} u_t
     + 2 \phi_{t} \phi_{tu} + \phi_{xtt})
\vspace{1mm}\\
\quad
- \alpha \phi (\phi_{xx} + \phi_{u} \phi_{x} 
+ \phi \phi_{u}^{2} + \phi^{2} \phi_{uu}
     + 2 \phi \phi_{xu})
                  - \beta (\phi_{x} + \phi \phi_{u})^{2}
=0
\ea \label{delzstar}
\ee
which involves the inf\/initesimal $\phi$ that is to be determined.
As in the $\tau\not=0$ case we apply the classical Lie algorithm to
this equation using the second prolongation $ {\rm pr}^{(2)} {\bf v}$
and eliminate $u_{tt}$ using (\ref{delzstar}). Similar to the
nonclassical method in the generic case $\tau\not=0$, when $\mu\not=0$
the coef\/f\/icient of the highest  derivative term, $u_{tt}$ is not
necessarily zero, thus singular cases are induced. As in the previous
section we consider the cases (i) $\mu=0$ and (ii) $\mu\not=0$ separately.

\subsection{Case (i) $\mu=0$.}

Generating the \deqs, again using {\tt symmgrp.max}, yields three equations,
the f\/irst
two being $\phi_{uu}=0$, $\phi_{tu}=0$. Hence we look for solutions like
$\phi=A(x) u+B(x,t)$ in the third. Taking coef\/f\/icients
of powers of $u$ to be zero yields a system of
three equations in $A$ and $B$.
\be
\ba{l}
\alpha A A_{xxx}
+ 2 \beta A^{2} A_{xx} + A_{xxxx} + \alpha A_{x} A_{xx}
+ 5 \beta A A_{x}^{2} + 6 \alpha A A_{x}^{2} + 10 A_{xx} A_{x}
\vspace{1mm}\\
\quad + 2 \gamma A_{xx} + 5 A A_{xxx} + \alpha A^{5} + 10 A^{2} A_{xx}
      + 10 A^{3} A_{x} + A^{5} + \beta A^{5}+ 15 A A_{x}^{2}
\vspace{1mm}\\
\quad + 4 \gamma A^{3}+ 2 \beta A_{x} A_{xx} + 10 \gamma A A_{x}
+ 4 \alpha A^{2} A_{xx}
+ 6 \beta A^{3} A_{x} + 7 \alpha A^{3} A_{x}=0,
\ea\label{abeqni}
\ee
\be
\ba{l}
 5 \alpha B A_{x}^{2}
+ 2 \beta A^{2} B_{xx} + \alpha B A_{xxx} + 2 \beta B A^{4}
     + 13 B A A_{xx} + 2 \alpha A^{3} B_{x} + 10 \gamma B A_{x}
\vspace{1mm}\\
\quad + 7 A A_{x} B_{x} + 2 \beta A_{x} B_{xx} + \alpha A_{x} B_{xx}
     + 2 \kappa A A_{x} + \alpha A^{2} B_{xx} + 15 B A^{2} A_{x}
\vspace{1mm}\\
\quad  + \alpha B_{x} A_{xx} + 6 \gamma A B_{x}
+ 8 \gamma B A^{2} + 2 \beta A^{3} B_{x}
   + 2 \beta B_{x} A_{xx} + 2 \alpha B A^{4}
 + \alpha A B_{xxx}
\vspace{1mm} \\
 \quad + B_{xxxx}  + 6 A_{xx} B_{x} + 2 \gamma B_{xx} + 5 B A_{xxx}
 + A B_{xxx} + \kappa A_{xx} + 2 B A^{4} + A^{2} B_{xx}
\vspace{1mm} \\
 \quad  + A^{3} B_{x}
     + 4 A_{x} B_{xx} + 11 B A_{x}^{2} + 4 \beta B A_{x}^{2}
     + 6 \beta A B_{x} A_{x} + 7 \alpha A B_{x} A_{x} + 7 \alpha B A A_{xx}
\vspace{1mm}\\
\quad + 2 \beta B A A_{xx} + 10 \beta B A^{2} A_{x} + 12 \alpha B A^{2} A_{x}
=0,
\ea \label{abeqnii}
\ee
\be
\ba{l}
\beta A B_{x}^{2} + 4 \gamma B^{2} A + 2 \beta B_{x} B_{xx} + 5 B^{2} A A_{x}
     + \alpha A B_{x}^{2} + 6 \gamma B B_{x} + 2 \kappa B A_{x} 
\vspace{1mm}\\
\quad + B A B_{xx}
+ \alpha B^{2} A^{3} + \beta B^{2} A^{3} + \alpha B_{x} B_{xx}
     + 3 \alpha B^{2} A_{xx} + \alpha B B_{xxx} + B A^{2} B_{x}
     \vspace{1mm}\\
\quad + 3 B A_{x} B_{x}
- B_{tt} + B B_{xxx} + \kappa B_{xx} + B^{2} A^{3}
     + 3 B^{2} A_{xx} + 2 \beta B A B_{xx} 
\vspace{1mm}\\
\quad + \alpha B A B_{xx}
 + 4 \beta B A_{x} B_{x}
+ 5 \alpha B A_{x} B_{x} + 2 \beta B A^{2} B_{x}
 + 2 \alpha B A^{2} B_{x}
 \vspace{1mm}\\
 \quad + 4 \beta B^{2} A A_{x} + 5 \alpha B^{2} A A_{x}
=0,
\ea \label{abeqniii}
\ee
We try to solve this system using the {\tt diffgrob2} package
interactively, however the expression swell is too
great to obtain meaningful output. Thus we proceed by making
ans\"atze on the form of $A(x)$, solve (\ref{abeqniii}) (a linear
equation in $B(x,t)$) then f\/inally (\ref{abeqnii}) gives the full picture.
Many solutions have been found as (\ref{abeqni}) lends itself to many 
ans\"atze through choices of parameter values. We present some
in  \S~\sen.3.

{\advance\topsep-1.5pt

\subsection{Case (ii) $\mu\not=0$.}
 
The nonclassical method, when the coef\/f\/icient of $u_{tt}$ is non-zero,
generates a system of three \deqs. However, far from being single-term
equations the f\/irst two contain 41 and 57 terms respectively,
and the third 329. The intractability of f\/inding all
solutions is obvious. To f\/ind some, we return to our previous case and
look for $\phi=A(x) u+B(x,t)$. Three equations then remain,
similar to (\ref{abeqni},\ref{abeqnii},\ref{abeqniii}) which we
tackle in the same vein as previously.
Some solutions are presented in \S~\sen.3.

As mentioned in the start of this section, singular solutions may exist, when the
coef\/f\/icient of $u_{tt}$ equals zero, i.e.\ when
\[
1 - \phi \phi_{uu} -\phi_u^2 - \phi_{ux} =0. \label{sing}
\]
This may be integrated \wrt\ $u$ to give
\[
u- \phi\phi_u - \phi_x = H(x,t). \label{singint}
\]
If $\phi$ satisf\/ies (\ref{sing}) then the coef\/f\/icients of
$u_t^2$ and $u_t$ in (\ref{delzstar}) are both zero. Since
no $u$-derivatives now exist in (\ref{delzstar}) what is left
must also be zero, i.e.
\[
\ba{l}
(2\gamma +\alpha) \phi^2 -\alpha H_x \phi + (2\gamma +
\beta+1) u^2
\vspace{1mm}\\
\qquad +\left[\kappa - (2\gamma+2\beta+1) H -H_{xx}\right]u
+\beta H^2 -\kappa H - H_{tt} =0.
\ea \label{kone}
\]
Thus we need to solve (\ref{singint}) and (\ref{kone}). Note that
once we have found $\phi(x,t,u)$, the related exact solution  
is given by solving the \isc , with no further restrictions
on the solution. The following are distinct from
each other and from solutions in  \S~\sen.3.

\smallskip
\noindent
{\bfseries \itshape Case (a)} $\gamma=-\tfr12$, $\alpha=1$ and $\beta=\kappa=0$.
In this case $\phi(x,t,u)$ is given by the relation
\[
u- \phi \phi_u - \phi_x = c_1 t + c_2.
\]
For instance, if $\phi (x,t,u)$ is linear in $u$
we have the exact solution
\[
u(x,t) = w(t) \cosh [x+ A(t)] + B(t) \sinh [x+ A(t)]  +c_1 t+c_2,
\]
where $w(t)$, $A(t)$ and $B(t)$ are \afs.

\smallskip
\noindent
{\bfseries \itshape Case (b)} $\alpha=-2 \gamma$, $\beta = -1-2 \gamma$ and $\gamma\not=-\tfr12$.
In this case $\phi (x,t,u)$ is given by the relation
\[
u- \phi \phi_u - \phi_x = - {\kappa \over 1+2 \gamma} .
\]
For instance, if $\phi (x,t,u)$ is linear in $u$ we have the exact solution
\[
u(x,t) = w(t) \cosh [x+ A(t)] + B(t) \sinh [x+ A(t)] -
{\kappa \over 1+2 \gamma},
\]
where $w(t)$, $A(t)$ and $B(t)$ are \afs.

\smallskip
\noindent
{\bfseries \itshape Case (c)} $\alpha=-2 \gamma$ and $\beta = -1-2 \gamma$. In this case
\[
\phi(x,t,u)= \pm u -\kappa x \mp
\left(- \tfr12 \kappa t^2 +c_1 t+c_2\right)\pm \kappa,
\]
and so 
\[
u(x,t) = w(t) \exp (\pm x) \pm \kappa x+ \left(- \tfr12 \kappa t^2 +c_1 t+c_2
\right),
\]
where $w(t)$ is an \af.

\smallskip
\noindent
{\bfseries \itshape Case (d)}
$\gamma=-\tfr12$, $\alpha=1$ and $\beta=0$. In this case
\[
\phi(x,t,u) = {\kappa u - H_{xx} u - \kappa H -H_{tt} \over H_x},
\]
where $H_x (x,t) \not=0$ and also $H(x,t)$ satisf\/ies the system
\[
\kappa H_{xx} +H_x^2 -\kappa^2 =0 \qquad
\left(\kappa^2-H_x^2\right) H_{tt} -2 \kappa H_{xt}^2
+\kappa^2\left(\kappa^2-H_x^2\right) =0 .
\]
We have assumed that $\kappa^2-H_x^2\not=0$,
for a dif\/ferent solution to (c). This yields
\[
\ba{l}
u(x,t) = \left[w(t) - 2 \kappa x\right] \sinh z\cosh  z
\vspace{2mm}\\
\qquad + \cosh^2 z
\left[ 4 \kappa \log ( \cosh  z ) - \kappa^2 t^2 +2 c_3 t +2 c_4 -
2\kappa +2 c_1^2\right] - c_1^2,
\ea
\]
where $z= \tfr12 (x+c_1 t +c_2)$ and $w(t)$ is an \af.

}   

\smallskip\noindent
{\bfseries \itshape Case (e)}
$\beta =-1-\alpha -4 \gamma$, $\gamma\not=\tfr12$ and $\alpha+ 2 \gamma\not=0$. In this case
\[
\phi(x,t,u)= \pm \left( u+ {\kappa \over 1+2 \gamma} \right),
\]
and so
\[
u(x,t) = w(t) \exp(\pm x) - {\kappa \over 1+2 \gamma},
\]
where $w(t)$ is an \af. 

\smallskip\noindent
{\bfseries \itshape Case (f)} $\gamma=0$, $\alpha=-2$ and $\beta=1$.
In this case $\phi(x,t,u)$ satisf\/ies
\[
-2 \phi^2 +2 \phi H_x + 2 u^2 +2 \kappa u - 2 H u +\kappa^2 + H^2 =0
\]
and $H(x,t)$ satisf\/ies the system
\[
H_{xx} + \kappa + H =0, \qquad H_{tt} + \kappa^2 +\kappa H =0.
\]
Then 
\[
u(x,t) = \tfr12 (A^2 + B^2)^{1/2}\sinh
[\pm x + w(t) ] - \kappa + \tfr12 (A \sin x + B \cos x),
\]
where $A(t)$ and $B(t)$ satisfy
\[
\frac{\d^2 A}{\d t^2} + \kappa A=0,
\qquad \frac{\d^2 B}{\d t^2} + \kappa B =0
\]
and $w(t)$ is an \af.

\smallskip\noindent
{\bfseries \itshape Case (g)}
$\gamma=0$, $\beta =-1-\alpha$, $\kappa =0$ and $\alpha=(c_1-1)^2
/c_1$, where $c_1\not=0,1 $. In this case
\[
\phi(x,t,u) = u + {c_2 t+c_3 \over c_1 -1} \exp(-c_1 x),
\]
and so 
\[
u(x,t) =\left\{
\begin{array}{ll}
w_1(t) \e^x - \tfr12 (c_2 t+c_3)x \e^x,
& \mbox{if} \quad c_1=-1,
\vspace{2mm}\\
  w_2(t) \e^x +\left(1-c_1^2\right)^{-1} (c_2 t+c_3)\e^{-c_1 x},&
  \mbox{if} \quad c_1\not=-1,
\end{array}\right.
\]
where $w_1(t)$ and $w_2(t)$ are \afs.

\smallskip\noindent
{\bfseries \itshape Case (h)}
$\gamma=0$, $\beta =-1-\alpha$, $c_1^2 +\alpha^2 c_1 +2 c_1 +4 \alpha
 c_1 +1=0$, $\alpha\not=0,-2$ and $c_1\not=0,1$.
 In this case $\phi(x,t,u)$ satisf\/ies
\[
\alpha \phi^2 - \alpha H_x \phi -  \alpha u^2 +u\left[\kappa
+H(1+2 \alpha) - H_{xx}\right] - \kappa H - (\alpha +1) H^2
 - H_{tt} =0,
\]
where $H(x,t)$ satisf\/ies the system 
\[
H_{tt} + c_1 \kappa H =0 \qquad (\alpha+2) H_x \pm (1-c_1 )(H+\kappa) =0.
\]
Then
\[
u(x,t) = \left[ {1+2\alpha +c_1 \over 2\alpha} w(t) +
 g(x,t) \right] \exp \left\{ {(c_1 -1) x \over \alpha+2}\right\} - \kappa,
\]
where $w(t)$ satisf\/ies
\[
\wtt+c_1 \kappa w=0,
\]
and
$g(x,t)$ satisf\/ies
\[
\ba{l}
\ds \alpha g_x^2 + {  2\alpha (c_1-1) \over \alpha+2} gg_x -
\alpha (c_1+1) g^2 +  (\alpha c_1 +c_1+1 )
 w(t) g_x
 \vspace{3mm}\\
\ds  \qquad - c_1 (1+\alpha +c_1) w(t) g + {c_1 (c_1-1) (\alpha+2)
\over 4 \alpha} w^2 (t) =0 .
\ea
\]

\subsection{Exact solutions} 
In this subsection some exact solutions are presented.
The \inl\ $\phi(x,t,u)$ is given, possibly up to
satisfying some equations, and then the solution, found
by solving the \isc\ (\ref{invsc}). 

\subsubsection{\mbox{$\ds \gamma=0$ and
$\ds\phi={u \over x} +H_1 (t) x + 3H_2 (t) x^3+
{H_3 (t) \over x} + H_4 (t) x^{2-\alpha}$.}}

Solving the \isc\ gives
\[
 u(x,t) = \left\{
 \ba{ll}
 \displaystyle xw(t) +H_1(t) x^2 + H_2(t) x^4 -H_3(t)
+ {H_4(t) x^{3-\alpha}\over 2-\alpha}, &
\mbox{if} \quad \alpha\not=2,
\vspace{3mm}\\
\ds x \~w (t) +H_1(t) x^2 + H_2(t) x^4 -H_3(t) +
H_4(t) x \log x, & \mbox{if} \quad \alpha=2.
\ea \right.
\]
Various types of solution are found, as seen in Table~4.1. The
$H_i(t)$ are obtained from by the \deqs, $w(t)$ by substituting
back into (\ref{fulleqn}).

\subsubsection{\mbox{$\phi=B(x,t)$.}}

\smallskip\noindent
{\bfseries \itshape Case (a)} $\gamma=0$. In this case
$B(x,t)= 4H_1 (t) x^3 +3H_2 (t) x^2 +2H_3(t) x +H_4 (t)$, where
$H_1 (t)$, $H_2 (t)$, $H_3(t)$ and $H_4 (t)$ satisfy
\[
\ba{l}
\ds \Htt{1}-24(6\beta +4\alpha +1) H_1^2 =0,
\vspace{3mm}\\
\ds  \Htt{2}-24(6\beta +4\alpha +1) H_1 H_2 =0,
\vspace{3mm}\\
\ds  \Htt{3}-24(2\beta +2\alpha +1) H_1 H_3 =18 (2\beta +\alpha) H_2^2
 +12 \kappa H_1 +12 \mu \Htt{1},
 \vspace{3mm} \\
\ds \Htt{4} -24(\alpha +1) H_1 H_4=12(2\beta +\alpha) H_2 H_3 +
6\kappa H_2 +6 \mu \Htt{2}.
\ea
\]
Then
\[
u(x,t) = w(t) +H_1(t) x^4 +H_2(t) x^3+ H_3(t) x^2 +H_4(t) x,
\]
where $w(t)$ satisf\/ies
\[
 \wtt-24 H_1 w =2 \kappa H_3 +6 \alpha H_2 H_4 +4\beta H_3^2
+ 2 \mu \Htt{3}.
\] 

\begin{center}
{\bf Table 4.1}

\medskip
\begin{tabular}{|c|c|}
\hline & \\[-5pt]
Parameters & $H_i(t)$ and $w(t)$  satisfy \\[5pt] \hline & \\[-5pt]
$\alpha=\tfr12$, $\beta=-\tfr18$&
$\ba{c}\displaystyle  H_1 =4 \kappa,\qquad H_2=H_3=\Htt{4}=0\\[5pt]
\displaystyle 32 \wtt = -5 H_4^2 \ea$\\
& \\[-5pt] \hline & \\[-5pt]
$\alpha=\tfr12$, $\kappa=0$ &
$\ba{c}\displaystyle  H_1 = \mu H_2 = H_3 =0 \vspace{3mm}\\
\displaystyle \Htt{2}-72 (1+2 \beta) H_2^2 =0 \vspace{3mm}\\
\displaystyle 16 \Htt{4} -3(480 \beta +303) H_2 H_4 =0 \vspace{3mm}\\
\displaystyle 8 \wtt -288 H_2 w = (50 \beta+5) H_4^2 \ea$\\
&   \\[-5pt] \hline & \\[-5pt]
$\ba{c} \ds\beta=\frac{\alpha^2-\alpha}{3-\alpha}\\[5pt]
\kappa =0,\quad\alpha\not=3 \ea$ &
$\ba{c}  H_1 = \mu H_2 = H_3 =0 \vspace{3mm}\\
\displaystyle (\alpha-3) \Htt{2}+24(2 \alpha+3)(\alpha+1) H_2^2=0 \vspace{3mm}\\
\displaystyle \Htt{4}+3(\alpha+2)(\alpha+1)(\alpha^2-\alpha-4) H_2 H_4 =0
\vspace{3mm}\\
\displaystyle \wtt - 24(\alpha+1) H_2 w=0 \vspace{3mm}\\
\displaystyle \wwtt- 72 H_2 \~w = 90 H_2 H_4
(\mbox{if}\enskip \alpha=2)  \ea$ \\
& \\[-5pt] \hline & \\[-5pt]
$\alpha=-2,\quad \beta =\tfr65$ &
$\ba{c} \displaystyle H_1 = \tfr56\kappa, \qquad \mu H_2 =\Htt{4}=0
\vspace{3mm}\\
\displaystyle 5 \Htt{2}- 24H_2^2 =0 \vspace{3mm}\\
\displaystyle 5 \Htt{3}-120 H_2 H_3 = -25 \kappa^2 \vspace{3mm}\\
\displaystyle \wtt +24 H_2 w = -30 H_3 H_4  \ea$ \\
& \\ \hline 
\end{tabular}
\end{center}

\smallskip\noindent
{\bfseries \itshape Case (b).} In this case
$B(x,t)= H_1 (t) +2H_2(t) x$ where
$H_1 (t)$ and $H_2 (t)$ satisfy
\be
\Htt{2} - 12\gamma H_2^2 =0, \label{Blini}
\ee
\be
\Htt{1} - 12 \gamma H_2 H_1 =0. \label{Blinii}
\ee
Then 
\[
u(x,t) =  w(t) + H_1(t) x +H_2(t) x^2,
\]
where $w(t)$ satisf\/ies
\be
 \wtt - 4\gamma H_2 w = 2\kappa H_2 +4(6 \gamma \mu +\beta)  H_2^2 +
2 \gamma H_1^2 .\label{Bliniii}
\ee

\smallskip\noindent
{\bfseries \itshape Case (c)} $\beta=1-\alpha$.
In this case $B(x,t)= c H_1 (t) \e^{c x} +c H_2 (t) \e^{-c x}
+ H_4 (t)$, with $c^2= -2 \gamma$ and where
$H_1 (t)$, $H_2 (t)$, $H_3(t)$ and $H_4 (t)$ satisfy
\be
\Htt{4}=0, \label{Bexpi}
\ee
\be
 (1+2\gamma \mu) \Htt{1}- 2\gamma c(2-\alpha)
H_4 H_1 + 2 \gamma \kappa H_1=0, \label{Bexpii}
\ee
\be
 (1+2\gamma \mu) \Htt{2}+2\gamma c(2-\alpha)
H_4 H_2 + 2 \gamma \kappa H_2=0. \label{Bexpiii}
\ee
Then
\[
u(x,t) = w(t) +H_1 (t) \e^{c x} -H_2 (t) \e^{-c x} +H_4(t) x,
\]
where $w(t)$ satisf\/ies
\beq
\wtt = 2 \gamma H_4^2 -16 \gamma^2 (1-\alpha) H_1 H_2 .\label{Bexpiv}
\eeq

\smallskip\noindent
{\bfseries \itshape Case (d)} $\alpha=2$ and $\beta=-1$.
In this case $B(x,t)=c H_1 (t) \e^{c x} +c H_2 (t) \e^{-c x}
+2 H_3 (t)x+ H_4 (t)$, with $c^2= -2 \gamma$ and where
$H_1 (t)$, $H_2 (t)$, $H_3(t)$ and $H_4 (t)$ satisfy
\be
\Htt{3} -12 \gamma H_3^2 =0, \label{Bexpli}
\ee
\[
\Htt{4} -12 \gamma H_3 H_4 =0, \label{Bexplii}
\]
\[
(1+2\gamma \mu) \Htt{1}-12 \gamma H_3 H_1 + 2 \kappa \gamma H_1=0,
\label{Bexpliii}
\]
\be
(1+2\gamma \mu) \Htt{2}-12 \gamma H_3 H_2 + 2 \kappa \gamma H_2=0.
\label{Bexpliv} 
\ee
Then
\[
u(x,t) = w(t) +H_1(t) \e^{c x} -H_2(t) \e^{-c x} + H_3(t) x^2 +H_4(t) x,
\]
where $w(t)$ satisf\/ies
\beq
\wtt -4\gamma H_3 w= 2\kappa H_3 +2\gamma H_4^2 +4(6 \gamma \mu-1)
 H_3^2 +16\gamma^2 H_1 H_2.\label{Bexplv}
 \eeq

\subsubsection{\mbox{$\gamma= \beta+\alpha+1=0$ and
$\phi= R (u +H_1 (t) + H_2(t) \e^{Rx} + H_3
(t) \e^{m_+ x} + H_4(t) \e^{m_- x})$.}}

Here $m_{\pm} = -\tfr12 R (2+\alpha\pm n)$, $n=\sqrt{\alpha
(\alpha+4)}$, with $R\not=\pm1$ a non-zero constant. Solving the \isc\ yields
\[
\ba{l}
\ds u(x,t) =
 w(t) \e^{Rx} - H_1(t)+R H_2(t) \e^{Rx}
 \vspace{3mm}\\
\ds  \qquad - {2 H_3(t) \over
 4+\alpha + n} \exp \left\{ -\tfr12 R x (2+\alpha + n)
\right\} - {2 H_4 (t) \over
 4+\alpha - n} \exp \left\{ -\tfr12 R x (2+\alpha - n)
\right\}.
\ea
\]
The solutions are represented in Table~4.2

The equations that the various
$H_i(t)$ satisfy in this subsection are all solvable,
and the order in which a list of equations should be solved is from
the top down. The only nonlinear equations all have either polynomial
solutions (sometimes only in special cases of the parameters) or are
equivalent to the Weierstrass elliptic function equation
(\ref{ellip}). The homogeneous part of any linear equation is either 
of Euler-type, is equivalent to the Airy equation \cite{refAS},
\[
 \Htt{}(t) + t H(t) =0 
\]
or is equivalent to the Lam\'e equation \cite{refI}
\beq 
\Htt{}(t) -\{ k +n(n+1) \wp (t)\} H(t) =0 .\label{lame}
\eeq
The particular integral of any non-homogeneous linear equation may
always be found, up to quadratures, using the method of variation of
parameters. 

For instance consider the solution of {\sl 4.3.2} case (b) above.
There are essentially two separate cases to consider, either (i)
$\gamma=0$ or (ii) $\gamma\not=0$. 

\smallskip\noindent 
{\bfseries \itshape Case (i)} $\gamma=0$. The functions $H_1(t)$ and
$H_2(t)$ are trivially found from (\ref{Blini}) and (\ref{Blinii}) to
be $H_1(t)=c_1t+c_2$ and $H_2(t)=c_3t+c_4$, 
then (\ref{Bliniii}) becomes
\[
\wtt = 2\kappa (c_3t+c_4) + 4 \beta (c_3t+c_4)^2
\]
which may be integrated twice to yield the exact solution
\[
 u(x,t) = \left\{ 
\ba{ll}
\ba{l}
\displaystyle {\kappa \over 3c_3^2}
 (c_3t+c_4)^3 + {\beta  \over 3c_3^2} (c_3t+c_4)^4
\vspace{3mm}\\
\ds \hspace*{3cm}  + c_5 t+c_6 +(c_1t+c_2)x+(c_3t+c_4) x^2,
\ea  & \mbox{if} \quad c_3\not=0,
\vspace{3mm}\\
\ds  \left(\kappa c_4 + 2 \beta c_4^2 \right) t^2 
+ c_5 t+c_6 +(c_1t+c_2)x+c_4 x^2,
 & \mbox{if} \quad  c_3=0. 
\ea\right.
\]

\smallskip\noindent 
{\bfseries \itshape Case (ii)} $\gamma\not=0$. Equation (\ref{Blini}) 
may be transformed into the Weierstrass elliptic function
equation (\ref{ellip}), hence $H_2(t)$ has solution
$H_2(t)=\wp (t+t_0;0,g_3) /(2\gamma)$. Now $H_1(t)$
satisf\/ies the Lam\'e equation 
\[
\Htt{1} - 6 \wp (t+t_0;0,g_3) H_1 =0,
\]
which has general solution
\[
H_1(t) = c_1 \wp (t+t_0;0,g_3) + c_2 \wp (t+t_0;0,g_3)
 \int^{t+t_0} {\d s \over \wp^2 (s;0,g_3)},
\]
where $c_1$ and $c_2$ are \acs. Now $w(t)$ satisf\/ies
the inhomogeneous Lam\'e equation
\beq
\wtt -2 \wp (t+t_0;0,g_3) w = Q(t),
 \label{wheqn}
\eeq
where $Q(t)=2\kappa H_2(t) +4(6 \gamma \mu +\beta)  H_2^2(t) +
2 \gamma H_1^2(t)$, with $H_1(t)$ and $H_2(t)$ as above. The 
general solution of the homogeneous part of 
this Lam\'e equation is given by
\[
w_{\rm CF} (t) = c_3 w_1 (t+t_0) + c_4 w_2 (t+t_0),
\]
where $c_3$ and $c_4$ are \acs,
\[
 w_1 (t) = \exp \{ -t \zeta (a) \} {\sigma (t+a) \over
 \sigma (t)}, \qquad w_2(t) = \exp \{ t \zeta (a) \}
{\sigma (t-a) \over \sigma (t)}
\]
in which $\zeta (z)$ and $\sigma (z)$ are the Weierstrass
zeta and sigma functions def\/ined by the dif\/ferential
equations
\[
\frac{\d \zeta}{\d z} = -\wp (z), \qquad
 \frac{\d}{\d z} \log \sigma (z) = \zeta (z)
\]
together with the conditions
\[
 \lim_{z\to0} \left( \zeta(z) - \frac1z \right) =0,
 \qquad \lim_{z\to0} \left( {\sigma(z)\over z} \right) =1
\]
respectively (cf.\ \cite{refWW}), and $a$ is any solution of the 
transcendental equation
$$ \wp (a) =0 $$
i.e., $a$ is a zero of the Weierstrass elliptic function
(cf.\ \cite{refI}, p.379). Hence the general
solution of (\ref{wheqn}) is given by
\be
\ba{l}
\ds  w(t) = c_3 w_1 (t+t_0) + c_4 w_2 (t+t_0)
\vspace{3mm}\\
\ds \qquad+ \frac1{W(a)} \int^{t+t_0} \left[ w_1 (s) w_2 (t+t_0) -
w_1 (t+t_0)  w_2 (s)\right] Q(s)\, \d s,
\ea
 \label{wsolna} 
\ee
where $W(a)$ is the non-zero Wronskian
\beq 
W(a) = w_1 w_2' - w_1' w_2 = - \sigma^2 (a) \wp ' (a) 
\label{wsolnb}
\eeq
and $Q(t)$ is def\/ined above. We remark that in order to verify that 
(\ref{wsolna},\ref{wsolnb}) is a solution of (\ref{wheqn}) one uses
the following addition theorems for Weierstrass elliptic, zeta and
sigma functions 
\[
\ba{l}
\ds \zeta (s\pm t) = \zeta (s) \pm \zeta (t) +
\frac12 \left[ {\wp' (s) \mp \wp' (t) \over \wp (s) - \wp (t)}
\right], 
\vspace{3mm}\\
\ds \sigma (s+t) \sigma (s-t) = -\sigma^2 (s) \sigma^2 (t) [ \wp (s)
- \wp (t) ] 
\ea
\]
(cf.\ \cite{refWW}, p.451).

\begin{center}
{\bf Table 4.2}

\medskip

\begin{tabular}{|l|c|}
\hline & \\[-5pt]
Parameters& $H_i (t)$ and $w(t)$ satisfy  \\[5pt] \hline & \\[-5pt]
$\ba{l} \alpha=-4\\ \beta=3 \ea $ &
$\ba{c}\ds H_3 = H_4=\Htt{1} =0\vspace{3mm}\\
\ds (1-\mu R^2) \Htt{2} +R^2 (R^2 H_1 -\kappa) H_2 =0\vspace{3mm}\\
\ds (1-\mu R^2) \wtt+R^2 (R^2 H_1 -\kappa) w
 =2 \kappa R^2 H_2 - 4 R^4 H_1 H_2
+2 \mu R^2 \Htt{2} \ea$\\
& \\[-5pt] \hline & \\[-5pt]
$\ba{l} \alpha\enskip \mbox{arbitrary} \\[5pt]
j=\tfr72 \pm \tfr12 \\[5pt]   i=\tfr72 \mp \tfr12
\ea$ & $\ba{c} \ds H_2 = H_j=\Htt{1} =0 \vspace{3mm}\\
\ds (4- \mu (2+\alpha \pm n)^2 R^2) \Htt{i}+R^2 H_i \vspace{3mm}\\
\ds \quad\times [  R^2 H_1 ((\alpha^2+4\alpha+2)(2+\alpha \pm n)^2 -4)
-\kappa (2+\alpha \pm n)^2]=0 \vspace{3mm}\\
\ds (1-\mu R^2) \wtt+R^2 (R^2 H_1 -\kappa) w=0 \ea $ \\
& \\[-5pt] \hline & \\[-5pt]
$\ba{l} \alpha=-3\\ \beta=2 \ea$ & $\ba{c}\ds  H_2 =\Htt{1} =0 \vspace{3mm}\\
\ds (2+\mu R^2 (1+\i \sqrt 3)) \Htt{3}-H_1 H_3 R^4 (1-\i \sqrt 3)+
\kappa H_3 R^2 (1+\i \sqrt 3)=0 \vspace{3mm}\\
\ds (2+\mu R^2 (1-\i \sqrt 3)) \Htt{4}-H_1 H_4 R^4
(1+\i \sqrt 3)+\kappa H_4 R^2 (1-\i \sqrt 3)=0 \vspace{3mm}\\
\ds (1-\mu R^2) \wtt+R^2 (R^2 H_1 -\kappa) w= 6 H_3 H_4 R^4
\ea$ \\
& \\[-5pt] \hline & \\[-5pt]
$\ba{l} \alpha=-1\\ \beta=0\ea $ & $\ba{c} \ds H_2 =\Htt{1} =0 \vspace{3mm}\\
\ds (2+\mu R^2 (1-\i \sqrt 3)) \Htt{3}-H_1 H_3 R^4 (1+\i \sqrt 3)
+\kappa H_3 R^2 (1-\i \sqrt 3)=0 \vspace{3mm}\\
\ds (2+\mu R^2 (1+\i \sqrt 3)) \Htt{4}-H_1 H_4 R^4
(1-\i \sqrt 3)+\kappa H_4 R^2 (1+\i \sqrt 3)=0 \vspace{3mm}\\
\ds (1-\mu R^2) \wtt+R^2 (R^2 H_1 -\kappa) w= 0 \ea$ \\
 & \\ \hline
\end{tabular}
\end{center}

\newpage

\section{Discussion}
This paper has seen a classif\/ication of symmetry reductions of the
nonlinear fourth order \pde\ (\ref{fulleqn}) using the classical Lie
method and the nonclassical method due to Bluman and Cole. The
presence of arbitrary parameters in (\ref{fulleqn}) has led to a
large variety of reductions using both symmetry methods for various
combinations of these parameters. The use of the MAPLE package {\tt
diffgrob2} was crucial in this classif\/ication procedure. In the
classical case it identif\/ied the special values of the parameters for
which additional symmetries might occur. In the generic nonclassical
case the f\/lexibility of {\tt diffgrob2} allowed the fully nonlinear
\deqs\ to be solved completely, whilst in the so-called $\tau=0$ case
it allowed the salvage of many reductions from a somewhat intractable
calculation. 

An interesting aspect of the results in this paper is that the class
of reductions given by the nonclassical method, which are not
obtainable using the classical Lie method, were much more plentiful
and richer than the analogous results for the generalized \ch\
equation (\ref{gch}) given in \cite{refCMP}. 

An interesting problem this paper throws open is whether
(\ref{fulleqn}) is integrable, or perhaps more realistically for
which values of the parameters is (\ref{fulleqn}) integrable. 
Ef\/fectively, in f\/inding the symmetry reductions of (\ref{fulleqn}),
we have provided a f\/irst step in using the Painlev\'e ODE test for
integrability due to Ablowitz, Ramani and Segur
\cite{refARS,refARSi}. However the presence of so many reductions
makes this a lengthy task and so the PDE test due to Weiss, Tabor and
Carnevale \cite{refWTC} is a more inviting prospect. It is likely
though that extensions of this test, namely ``weak Painlev\'e
analysis'' \cite{refRDG,refRRDG} and ``perturbative Painlev\'e
analysis'' \cite{refCFP} will be necessary (for instance see
\cite{refGP}). We shall not pursue this further here. 

The FFCH equation (\ref{ffch}) may be thought of as an integrable
generalization of the Korteweg-de Vries equation (\ref{kdv}).
Analogous integrable generalizations of the modif\/ied Korteweg-de
Vries equation 
\[
 q_t = q_x + q_{xxx} + 3\gamma q^2q_x,
\]
the nonlinear Schr\"odinger equation
\[
\i q_t + q_{xx} + |q|^2q=0
\]
and Sine-Gordon equation
\[
 q_{xt} = \sin q
\]
are given by
\be
u_t + \nu u_{xxt} = u_x + q_{xxx} + \gamma[(u+\nu u_{xx})(u^2+\nu
u_x^2)]_x,\label{eqgmkdv}
\ee
\be
\i u_t + \i u_x + \mu u_{xt} + u_{xx} + \kappa u|u|^2 - \i\kappa\mu
|u|^2u_x=0, \label{eqgnls}
\ee
\be
u_{xt} = \sin (u +\mu u_{xx})=0,\label{eqgsg}
\ee
respectively, where $\mu$ and $\kappa$ are arbitrary constants
\cite{refFokb,refFOR,refOR}. 

Recently Clarkson, Gordoa and Pickering \cite{refCGP} derived
$2+1$-dimensional generalization of the FFCH equation (\ref{ffch})
given by 
\beq
\tfr12u_{y}u_{xxxx}+u_{xy} u_{xxx}-\alpha 
\left( \tfr12u_{y} u_{xx} + u_{x} u_{xy} \right) 
+ u_{xxxt}-\alpha u_{xt} = 0,\label{ggch21} 
\eeq
where $\alpha$ is an arbitrary constant. The FFCH equation
(\ref{ffch}) and is obtained from (\ref{ggch21}) under the reduction
$\partial_y=\partial_x$, with $v=u_x$. The $2+1$-dimensional FFCH
equation (\ref{ffch}) has the non-isospectral Lax pair 
\[
\ba{l}
4\psi_{xx} = \left[ \alpha -\lambda 
\left(u_{xxx} -\alpha u_x  \right)\right]\psi, 
\vspace{1mm}\\
\psi_t = \lambda^{-1}\psi_y-\tfr{1}{2} u_y\psi_x +\tfr{1}{4} u_{xy} \psi,
\ea
\]
with $\lambda$ satisfying $\lambda_y=\lambda\lambda_t$. Clarkson,
Gordoa and Pickering \cite{refCGP} also derived a $2$-component
generalisation of the FFCH equation (\ref{ffch}) in $2+1$-dimensions
given by 
\beq
\ba{l}
u_{xxxt}-\alpha u_{xt} = -\tfr{1}{2}u_{y}u_{xxxx}-u_{xy}
u_{xxx}+\alpha \left(\tfr{1}{2} u_{y} u_{xx} + u_{x} u_{xy}\right)
-\kappa u_{xxxy} +v_{y}, 
\vspace{1mm}\\ 
 v_{t} = -v u_{xy} -\tfr{1}{2} v_{x} u_y,
\ea \label{ggch22}
\eeq
which has the Lax pair
\[
\ba{l}
4(1+\kappa\lambda)\psi_{xx} =
\left[\alpha-\lambda\left(u_{xxx}-\alpha u_x \right) 
-\lambda^2 v\right]\psi,
\vspace{1mm}\\
\psi_t = \lambda^{-1}\psi_y-\tfr{1}{2} u_y\psi_x +\tfr{1}{4} u_{xy} \psi,
\ea
\]
where the spectral parameter $\lambda$ satisf\/ies
$\lambda_y=\lambda\lambda_t$. 

We believe that a study of symmetry reductions of
(\ref{eqgmkdv},\ref{eqgnls},\ref{eqgsg},\ref{ggch21},\ref{ggch22})
would be interesting, though we shall not pursue this further here.

\section*{Acknowlegdements}
We thank Elizabeth Mansf\/ield for many interesting discussions. The
research of TJP was supported by an EPSRC Postgraduate Research
Studentship, which is gratefully acknowledged. 

\section*{Appendix A}
In this appendix we list the \deqs\ that are generated
in  \S~3.2 in the generic case when $\xi^2+u\not=0$.
\[
\ba{l}
\xi_u =0
\vspace{1mm}\\
\alpha \xi^{2} \phi_{u}
- 2 \alpha \xi \xi_{t} + 4 \phi_{uu} u^{2} - \alpha \phi
          + 4 \xi^{4} \phi_{uu} + \alpha \phi_{u} u + 8 \xi^{2} \phi_{uu} u
          - 2 \alpha \xi^{2} \xi_{x}=0
\vspace{1mm}\\
\beta \xi^{2} \phi_{u} + 3 \phi_{uu} u^{2} +
3 \xi^{4} \phi_{uu}
      + 6 \xi^{2} \phi_{uu} u
- \beta \phi + \beta \phi_{u} u - 2 \beta \xi^{2} \xi_{x}
      - 2 \beta \xi \xi_{t}=0
\vspace{1mm}\\
12 \xi^{2} \phi_{uuu} u + 3 \alpha \phi_{uu} u + 6 \phi_{uuu} u^{2}
 + 2 \beta \phi_{uu} u + 2 \beta \xi^{2} \phi_{uu} + 3 \alpha \xi^{2} \phi_{uu}
       + 6 \xi^{4} \phi_{uuu}=0
\vspace{1mm}\\
7 \xi^{2} \xi_{xx} u - 10 \xi \xi_{x}^{2} u - 2 \xi^{2} \xi_{t} \phi_{u}
 - 4 \xi \xi_{xt} u
         - \alpha \xi^{2} \phi_{x} - 2 \xi \phi \phi_{u}
         - 4 \xi_{t} \xi_{x} u - \alpha \phi_{x} u + 2 \xi_{t} \phi_{u} u
\vspace{1mm}\\
\quad     + 2 \xi^{3} \phi \phi_{uu}
 + 5 \xi \xi_{x} \phi + 2 \xi \phi_{tu} u
    - 6 \xi^{2} \phi_{xu} u + 4 \xi^{2} \xi_{t} \xi_{x} - 2 \xi \xi_{t}^{2}
 + 6 \xi_{xx} u^{2} + 2 \xi^{3} \phi_{tu}
+ \xi_{tt} u
\vspace{1mm}\\
\quad+ \xi^{4} \xi_{xx}
         + \xi^{2} \xi_{tt} + 4 \xi \xi_{x} \phi_{u} u - 4 \phi_{xu} u^{2}
         + 2 \xi \phi \phi_{uu} u - 2 \xi^{4} \phi_{xu} - 4 \xi^{3} \xi_{xt}
         - \xi_{t} \phi =0
\vspace{1mm}\\
\alpha \xi^{2} \phi_{uuu} + \alpha
\phi_{uuu} u + 2 \xi^{2} \phi_{uuuu} u
            + \xi^{4} \phi_{uuuu} + \phi_{uuuu} u^{2}=0
\vspace{1mm}\\
3 \alpha \xi_{xx} u- 4 \beta \xi^{2} \phi_{xu} - 12 \phi_{xuu} u^{2}
- 3 \alpha \phi_{xu} u + 6 \xi \phi_{tuu} u + 3 \alpha \xi^{2} \xi_{xx}
 - 18 \xi^{2} \phi_{xuu} u
 \vspace{1mm}\\
\quad + 9 \xi_{t} \phi_{uu} u - 4 \beta \phi_{xu} u
    - 6 \xi \phi \phi_{uu} + 2 \beta \xi_{xx} u + 6 \xi^{3} \phi \phi_{uuu}
      - 15 \xi^{3} \xi_{x} \phi_{uu} - 3 \xi^{2} \xi_{t} \phi_{uu}
\vspace{1mm}\\
\quad - 3 \alpha \xi^{2} \phi_{xu} + 12 \xi^{3} \phi_{u} \phi_{uu}
 - 3 \xi \xi_{x} \phi_{uu} u + 6 \xi \phi \phi_{uuu} u
      + 12 \xi \phi_{u} \phi_{uu} u
+ 2 \beta \xi^{2} \xi_{xx}
\vspace{1mm}\\
\quad - 6 \xi^{4} \phi_{xuu} + 6 \xi^{3} \phi_{tuu}=0
\ea\hspace{-3.33pt}
\]
\[
\ba{l}
2 \xi^{3} \phi_{xu} \phi_{xx} - \xi^{3} \xi_{xx}
\phi_{xx}     + 2 \xi \xi_{t} \phi \phi_{xxu} - 2 \xi \xi_{t} \phi_{t}
 + \xi_{xx} \phi_{xt} u + \xi^{3} \phi_{xxu} \phi_{x}
 + 2 \xi_{xxt} \phi_{x} u
 \vspace{1mm}\\
\quad + \phi^{2} \phi_{uu} u
  - 2 \xi_{t} \phi_{x} u
 + 2 \xi^{2} \phi \phi_{tu} - \xi_{xx} \phi \phi_{x}
 - \phi_{t} \phi_{xxu} u
 - 2 \xi_{x} \phi_{xxt} u
 - 2 \xi \xi_{t} \xi_{xx} \phi_{x}
 \vspace{1mm}\\ \quad
+ 2 \phi \phi_{tu} u
  - \kappa \phi_{xx} u
     - 2 \phi_{xt} \phi_{xu} u + 4 \xi_{xt} \phi_{xx} u
     - 2 \phi \phi_{xxtu} u + 4 \xi \xi_{t} \phi_{xu} \phi_{x}
+ 2 \xi \xi_{t} \phi_{xxt}
\vspace{1mm}\\
 \quad
 - 2 \xi^{2} \phi_{xt} \phi_{xu}
     + 4 \xi_{x} \phi_{t} u
       - 2 \xi^{2} \phi \phi_{xxtu}
     - 4 \xi^{2} \phi_{xtu} \phi_{x} + 4 \xi^{2} \xi_{xt} \phi_{xx}
     - \xi^{2} \phi^{2} \phi_{xxuu} \vspace{1mm}\\
  \quad
+ \phi \phi_{u} \phi_{xx}
     - \xi^{2} \kappa \phi_{xx} + \xi^{2} \xi_{xx} \phi_{xt}
 - \xi^{2} \phi_{t} \phi_{xxu} + \xi \phi \phi_{x} - 2 \phi \phi_{xu}^{2} u
 - 2 \xi^{2} \phi_{uu} \phi_{x}^{2}
 \vspace{1mm}\\
 \quad - 2 \xi_{x} \phi \phi_{xx} + 2 \xi_{t} \phi_{xxx} u
 + 2 \xi^{2} \xi_{x} \phi_{xu} \phi_{x}
     + \xi^{2} \phi^{2} \phi_{uu} - \phi^{2} \phi_{xxuu} u
     + 2 \phi \phi_{xu} \phi_{x}
 \vspace{1mm}\\
\quad
- \xi^{2} \phi_{xxxx} u  + 2 \xi^{2} \xi_{xxt} \phi_{x}
- 2 \xi^{3} \xi_{x} \phi_{x}
 - 2 \xi^{2} \phi \phi_{xu}^{2}
 - 2 \phi_{tu} \phi_{xx} u
     + 4 \xi_{x}^{2} \phi_{xx} u
- 2 \gamma \phi_{xx} u^{2}
\vspace{1mm}\\
\quad + 2 \xi^{2} \xi_{x} \phi_{t}
- 2 \xi^{2} \phi_{tu} \phi_{xx}
 - \phi \phi_{u} \phi_{xxu} u - \xi \phi \phi_{xxx}
 - 2 \xi \xi_{t} \phi \phi_{u}
 - 2 \phi_{uu} \phi_{x}^{2} u
     - 4 \phi_{xtu} \phi_{x} u
\vspace{1mm}\\
\quad  - 2 \phi \phi_{uu} \phi_{xx} u
 - 4 \phi \phi_{xuu} \phi_{x} u - \xi^{2} \phi_{xxtt} + \xi^{2} \phi_{tt}
     - \phi_{xxxx} u^{2} + 2 \xi \xi_{x} \phi_{xxx} u
     + \xi_{xx} \phi \phi_{xu} u
\vspace{1mm}\\
\quad  + \xi^{2} \xi_{xx} \phi_{u} \phi_{x}
     + \xi_{x} \xi_{xx} \phi_{x} u - \xi \xi_{xx} \phi_{xx} u
     + \xi \phi_{xxu} \phi_{x} u + \xi_{xx} \phi_{u} \phi_{x} u
- 2 \xi_{x} \phi \phi_{xxu} u
\vspace{1mm}\\
\quad
 - 2 \xi^{2} \phi \phi_{uu} \phi_{xx}
     + 4 \xi_{x} \phi \phi_{u} u + 2 \xi \phi_{xu} \phi_{xx} u
     - 2 \xi_{x} \phi_{u} \phi_{xx} u - 2 \xi^{2} \gamma \phi_{xx} u
- 2 \xi_{x} \phi_{xu} \phi_{x} u
\vspace{1mm}\\
\quad  - \xi^{2} \xi_{x} \xi_{xx} \phi_{x}
 + 2 \xi^{2} \xi_{x} \phi \phi_{u} + 2 \xi \xi_{t} \phi_{u} \phi_{xx}
     - 4 \xi \xi_{t} \xi_{x} \phi_{xx} - 4 \xi^{2} \phi \phi_{xuu} \phi_{x}
  - \phi_{xxtt} u
\vspace{1mm}\\
\quad  + \phi_{tt} u - \phi^{2} \phi_{u} + \phi \phi_{xxt}
   - \phi \phi_{t} + \phi^{2} \phi_{xxu}
     - 2 \xi^{2} \phi_{u} \phi_{xu} \phi_{x}
     - \xi^{2} \phi \phi_{u} \phi_{xxu}
- 4 \xi \xi_{x} \phi_{x} u
\vspace{1mm}\\
\quad
+ \xi^{2} \xi_{xx} \phi \phi_{xu} - 2 \phi_{u} \phi_{xu} \phi_{x} u=0
\vspace{1mm}\\
2 \xi^{3} \xi_{xx} \phi_{u}
- 4 \xi^{2} \xi_{xxx} u - 8 \xi \xi_{t} \xi_{x}^{2}
- 3 \xi^{3} \xi_{x} \xi_{xx}
+ 4 \xi \xi_{t} \xi_{xt} + \phi^{2} \phi_{uuu} u
     - \kappa \phi + 4 \xi_{xt} \xi_{x} u
\vspace{1mm}\\
\quad
+ 2 \phi_{tu} \phi_{u} u - 2 \xi \xi_{t} u
 - 4 \xi \phi_{xtu} u
- 4 \xi^{3} \phi \phi_{xuu} + 2 \beta \xi^{2} \phi_{xx}
     - 4 \xi^{2} \xi_{x} u + 4 \xi_{x} \gamma u^{2} + 2 \xi^{2} \gamma \phi
\vspace{1mm}\\
\quad
- 4 \xi_{xt} \phi_{u} u + 2 \phi \phi_{tuu} u - 4 \xi^{2} \xi_{x} \phi_{tu}
     - 4 \xi^{2} \xi_{xt} \phi_{u} - 5 \xi^{3} \phi_{uu} \phi_{x}
 - 8 \xi_{x}^{2} \phi_{u} u +
 2 \xi_{x} \phi_{u}^{2} u
 \vspace{1mm}\\
 \quad
- 2 \xi \xi_{t} \kappa
       - 4 \xi^{3} \phi_{u} \phi_{xu} + 7 \xi_{t} \xi_{xx} u
 + \xi^{2} \phi_{t} \phi_{uu} - 4 \xi \xi_{xx} \phi
+ \phi_{t} \phi_{uu} u + \xi^{2} \phi^{2} \phi_{uuu}
\vspace{1mm}\\
\quad
+ 2 \xi \xi_{xxt} u
- \xi^{2} \xi_{t} \xi_{xx}
+ 2 \xi^{2} \xi_{t} \phi_{xu} - 8 \xi_{t} \phi_{xu} u
     + 6 \xi^{3} \xi_{x} \phi_{xu} - 2 \xi_{x} \phi_{tu} u
- 2 \xi \xi_{t} \phi \phi_{uu}
\vspace{1mm}\\
\quad
+ 2 \xi^{2} \phi_{tu} \phi_{u}
 - 2 \xi \xi_{t} \phi_{u}^{2} + 7 \xi^{2} \phi_{xxu} u
     + 8 \xi \xi_{t} \xi_{x} \phi_{u} - 4 \xi \phi_{u} \phi_{xu} u
     - 2 \xi_{x} \phi \phi_{uu} u - 4 \xi \xi_{x} \phi_{xu} u
\vspace{1mm}\\
\quad + 2 \beta \phi_{xx} u
+ 2 \xi^{2} \phi \phi_{tuu} + 2 \xi \xi_{xx} \phi_{u} u
     + \xi^{2} \phi_{ttu} + 8 \xi_{x}^{3} u - 4 \xi^{3} \phi_{xtu}
     + 2 \xi^{3} \xi_{xxt} + 5 \xi \phi \phi_{xu}
\vspace{1mm}\\
\quad+ 2 \xi_{xt} \phi - 4 \xi \xi_{t} \gamma u
+ \xi^{4} \phi_{xxu} + 5 \xi \xi_{x} \xi_{xx} u
     - 5 \xi \phi_{uu} \phi_{x} u - 4 \xi \phi \phi_{xuu} u
 - 2 \xi \xi_{t} \phi_{tu} + \phi_{ttu} u
 \vspace{1mm}\\
 \quad - \phi^{2} \phi_{uu}
     + 3 \xi^{2} \phi \phi_{u} \phi_{uu} + 3 \phi \phi_{u} \phi_{uu} u
     - \phi \phi_{u}^{2} - \phi \phi_{tu} + 2 \xi_{x} \kappa u
     - 4 \xi^{2} \xi_{x} \phi \phi_{uu} + 4 \xi_{x} \phi \phi_{u}
\vspace{1mm}\\
\quad + 6 \phi_{xxu} u^{2} - 4 \xi_{xxx} u^{2} - 2 \xi^{2} \xi_{xtt}
     - 2 \xi^{4} \xi_{x} - 2 \xi_{xtt} u - 4 \xi_{x}^{2} \phi + \xi^{2} \phi
     + 8 \xi^{2} \xi_{xt} \xi_{x}=0
\vspace{1mm}\\
2 \xi \phi_{tuuu} u - 5 \xi^{3} \xi_{x} \phi_{uuu} - 3 \alpha \phi_{xuu} u
      + 2 \xi^{3} \phi \phi_{uuuu} + 6 \xi \phi_{uu}^{2} u
      + 3 \xi_{t} \phi_{uuu} u + 6 \xi^{3} \phi_{u} \phi_{uuu}
\vspace{1mm}\\
\quad  + 6 \xi^{3} \phi_{uu}^{2} + 2 \xi^{3} \phi_{tuuu}
- 2 \xi \phi \phi_{uuu}- \xi^{2} \xi_{t} \phi_{uuu}
- 3 \alpha \xi^{2} \phi_{xuu} + 6 \xi \phi_{u} \phi_{uuu} u
\vspace{1mm}\\
\quad + 2 \xi \phi \phi_{uuuu} u- 6 \xi^{2} \phi_{xuuu} u
- \xi \xi_{x} \phi_{uuu} u - 4 \phi_{xuuu} u^{2} - 2 \xi^{4} \phi_{xuuu}=0
\vspace{1mm}\\
 4 \xi^{2} \phi_{tu} \phi_{uu}
+ \xi^{2} \phi_{t} \phi_{uuu}
       - \alpha \xi^{2} \xi_{xxx} + 4 \xi^{2} \phi_{tuu} \phi_{u}
       - 6 \xi_{xt} \phi_{uu} u - 2 \gamma \phi + 4 \phi_{tu} \phi_{uu} u
\vspace{1mm}\\
\quad  + 2 \gamma \phi_{u} u
+ \phi^{2} \phi_{uuuu} u
+ 2 \phi \phi_{tuuu} u  - \alpha \xi_{xxx} u
+ 5 \xi \phi \phi_{xuu}+ 2 \xi^{2} \xi_{t} \phi_{xuu}
       + 4 \xi_{x} \phi \phi_{uu} \vspace{1mm}\\
\quad
- 4 \xi^{2} \xi_{x} \phi_{tuu}
+ 2 \xi^{2} \gamma \phi_{u} + 4 \xi^{2} \xi_{x}^{2} \phi_{uu}
       - 5 \xi^{3} \phi_{uuu} \phi_{x} + 2 \gamma \phi_{uu} u^{2}
       + \xi^{2} \kappa \phi_{uu} + 6 \xi^{3} \xi_{x} \phi_{xuu}
\vspace{1mm}\\
\quad - 8 \xi^{3} \phi_{u} \phi_{xuu}
 - 4 \xi^{3} \phi \phi_{xuuu}
 - 2 \xi \xi_{t} \phi_{tuu} + 4 \phi_{tuu} \phi_{u} u
       + \phi_{t} \phi_{uuu} u - 14 \xi^{3} \phi_{xu} \phi_{uu}
\vspace{1mm}\\
\quad  + 2 \xi^{2} \phi \phi_{tuuu}
- 8 \xi_{t} \phi_{xuu} u
       - 6 \xi^{2} \xi_{xt} \phi_{uu} + 6 \xi^{3} \xi_{xx} \phi_{uu}
 + 3 \alpha \phi_{xxu} u - 2 \xi_{x} \phi_{tuu} u  \vspace{1mm}\\
 \quad
- \xi^{2} \phi_{uu} u
       + 4 \phi \phi_{uu}^{2} u + \kappa \phi_{uu} u + 4 \xi_{x} \gamma u
       + 4 \phi_{u}^{2} \phi_{uu} u - 8 \xi \phi_{u} \phi_{xuu} u
 - 4 \xi \phi_{xtuu} u
 \vspace{1mm}\\
 \quad  + \xi^{2} \phi^{2} \phi_{uuuu}
  + 7 \xi^{2} \phi_{xxuu} u - 14 \xi \phi_{xu} \phi_{uu} u
       - 2 \xi_{x} \phi_{u} \phi_{uu} u + 2 \xi^{2} \gamma \phi_{uu} u
       + 6 \xi \xi_{xx} \phi_{uu} u
\vspace{1mm}\\
\quad
+ 4 \xi^{2} \phi_{u}^{2} \phi_{uu}
 - 4 \xi \xi_{t} \gamma - \phi^{2} \phi_{uuu}- \phi \phi_{tuu}
       + 6 \phi_{xxuu} u^{2} - 4 \xi^{3} \phi_{xtuu}
       - 2 \xi \xi_{t} \phi \phi_{uuu}
\vspace{1mm}\\
\quad  - 4 \xi \phi \phi_{xuuu} u
 - 4 \xi_{x}^{2} \phi_{uu} u - 5 \xi \phi_{uuu} \phi_{x} u
 - \xi^{4} \phi_{uu} - 8 \xi^{2} \xi_{x} \phi_{u} \phi_{uu}
       + 5 \phi \phi_{u} \phi_{uuu} u
\vspace{1mm}\\
\quad  + 5 \xi^{2} \phi \phi_{u} \phi_{uuu}
       + 8 \xi \xi_{t} \xi_{x} \phi_{uu} - 2 \xi_{x} \phi \phi_{uuu} u
       + 3 \alpha \xi^{2} \phi_{xxu} + \xi^{4} \phi_{xxuu}
 - 3 \phi \phi_{u} \phi_{uu}
 \vspace{1mm}\\
\quad  - 4 \xi \xi_{x} \phi_{xuu} u
  - 4 \xi^{2} \xi_{x} \phi \phi_{uuu} + \xi^{2} \phi_{ttuu} + \phi_{ttuu} u
+ 4 \xi^{2} \phi \phi_{uu}^{2} - 6 \xi \xi_{t} \phi_{u} \phi_{uu}=0 
\ea
\]

\newpage

\[
\ba{l}
 2 \phi_{xt} \phi_{uu} u + 4 \phi_{tuu} \phi_{x} u + 2 \xi_{t} \xi_{x} u
 + 4 \xi \phi \phi_{xxu} + 8 \xi \xi_{x} \phi_{u} u - 2 \xi_{xx} \gamma u^{2}
+ 4 \phi \phi_{xtuu} u
\vspace{1mm}\\
\quad - 6 \xi \phi_{xu}^{2} u
+ 6 \phi_{tu} \phi_{xu} u
     - 6 \xi^{3} \phi_{xuu} \phi_{x} + 2 \phi_{t} \phi_{xuu} u
     + \xi^{2} \xi_{t} \phi_{xxu} + 4 \gamma \phi_{x} u
\vspace{1mm}\\
\quad
+ 4 \phi_{xxxu} u^{2}
     + 2 \phi_{xttu} u
+ 6 \xi^{2} \phi \phi_{u} \phi_{xuu}
     - 4 \xi^{2} \xi_{x} \phi \phi_{xuu} + 4 \xi^{2} \phi_{xtu} \phi_{u}
     + 4 \xi^{2} \phi \phi_{xtuu} \vspace{1mm}\\
\quad
- 2 \phi^{2} \phi_{xuu}
     + 2 \xi_{xx} \phi \phi_{u} + 2 \xi^{2} \phi_{u}^{2} \phi_{xu}
     - 4 \xi^{3} \phi_{uu} \phi_{xx} + 2 \phi_{u}^{2} \phi_{xu} u
 + 5 \xi_{x}^{2} \xi_{xx} u - 2 \xi^{2} \xi_{t} \phi_{u}
 \vspace{1mm}\\
\quad
+ 2 \xi \phi_{tu} u
 + \xi_{tt} u + 2 \xi^{3} \phi \phi_{uu}
- 2 \xi^{3} \phi_{u} \phi_{xxu}
     - 4 \phi \phi_{u} \phi_{xu} - 10 \xi^{2} \xi_{xt} \phi_{xu}
     - 4 \xi \xi_{t} \phi_{xtu}
     \vspace{1mm}\\
\quad
 + 4 \xi^{2} \phi_{tuu} \phi_{x}
     + 2 \xi^{2} \xi_{t} \xi_{x}
+ 4 \gamma \phi_{xu} u^{2}
+ 5 \xi^{2} \xi_{xt} \xi_{xx}
     + \alpha \phi_{xxx} u + \xi \xi_{x} \phi
     + 2 \xi^{2} \phi^{2} \phi_{xuuu}
\vspace{1mm}\\
\quad
     - 4 \xi \xi_{t} \phi \phi_{xuu} 
+ 6 \xi^{2} \phi_{tu} \phi_{xu}
 + 5 \xi \xi_{xx} \phi_{xu} u
+ 2 \kappa \phi_{xu} u - 10 \xi_{xt} \phi_{xu} u
     + 4 \xi^{2} \gamma \phi_{x} - \xi_{xx} \kappa u
\vspace{1mm}\\
\quad
- 3 \xi_{x} \xi_{xx} \phi
     + 2 \xi^{2} \xi_{xxt} \xi_{x} + 2 \xi^{2} \kappa \phi_{xu}
     + \xi^{3} \xi_{x} \phi_{xxu} + 2 \xi^{2} \phi_{xt} \phi_{uu}
     + 2 \xi^{2} \xi_{x}^{2} \phi_{xu} - \xi \xi_{xxx} \phi
\vspace{1mm}\\
\quad
      - 2 \xi^{2} \xi_{xxt} \phi_{u} - 2 \phi \phi_{uu} \phi_{x}
     - \xi^{2} \xi_{xx} \phi_{u}^{2}
+ 2 \xi_{t} \phi_{u} u + 2 \xi_{t} \xi_{xxx} u
     + 5 \xi^{3} \xi_{xx} \phi_{xu} + 2 \xi \xi_{t} \xi_{xxt}
     \vspace{1mm}\\
\quad
- \xi^{2} \xi_{x}^{2} \xi_{xx}
  + 4 \phi \phi_{uuu} \phi_{x} u
     + 6 \phi_{u} \phi_{uu} \phi_{x} u - 2 \phi \phi_{xtu}
     + 8 \xi^{2} \phi \phi_{xu} \phi_{uu} - 2 \xi_{x} \xi_{xx} \phi_{u} u
\vspace{1mm}\\
\quad     - 3 \xi^{2} \xi_{xx} \phi_{tu}
- 2 \xi \phi_{xxtu} u - \xi \xi_{xx}^{2} u
- \xi_{xxtt} u  - 3 \xi^{2} \xi_{xx} \phi \phi_{uu}
 + 2 \xi^{3} \phi_{tu}  - \xi^{3} \xi_{xx}^{2}
 \vspace{1mm}\\
 \quad
- 2 \xi^{3} \phi \phi_{xxuu} - \xi^{2} \xi_{xxxx} u
     + 2 \xi^{2} \phi_{t} \phi_{xuu} - \xi_{xx} \phi_{u}^{2} u
      + 6 \xi_{x} \phi \phi_{xu} - 2 \xi_{xxt} \phi_{u} u
\vspace{1mm}\\
\quad  - 4 \xi^{2} \xi_{x} \phi_{xtu}
   - 2 \xi \phi \phi_{u} - \xi^{2} \xi_{xx} \kappa
     + 5 \xi_{xt} \xi_{xx} u + \alpha \xi^{2} \phi_{xxx}
     - 8 \xi \xi_{t} \phi_{u} \phi_{xu}- 4 \xi \xi_{x}^{2} u
\vspace{1mm}\\
\quad
  - 6 \xi^{3} \phi_{xu}^{2}
  + 6 \phi \phi_{u} \phi_{xuu} u
 + 8 \phi \phi_{xu} \phi_{uu} u - 4 \xi^{2} \xi_{x} \phi_{u} \phi_{xu}
 + \xi_{xxt} \phi - 2 \xi^{3} \phi_{xxtu}
 \vspace{1mm}\\
\quad  - 10 \xi_{x}^{2} \phi_{xu} u
  + 4 \phi_{xtu} \phi_{u} u - 2 \xi \xi_{t}^{2}
     - 4 \xi \xi_{t} \phi_{uu} \phi_{x} - 3 \xi_{xx} \phi_{tu} u
     + 2 \phi^{2} \phi_{xuuu} u
\vspace{1mm}\\
\quad  + 12 \xi \xi_{t} \xi_{x} \phi_{xu}
- \xi_{xxxx} u^{2} - 2 \xi_{x} \phi_{uu} \phi_{x} u - 2 \xi^{3} \xi_{x}^{2}
     + 2 \xi^{2} \phi_{xttu} + 4 \xi^{3} \xi_{x} \phi_{u}
     - 6 \xi \phi_{xuu} \phi_{x} u 
\vspace{1mm}\\
\quad
 + 2 \xi \phi \phi_{uu} u
     - 2 \xi \phi_{u} \phi_{xxu} u - 2 \xi \phi \phi_{xxuu} u
     + 4 \xi^{2} \gamma \phi_{xu} u + 2 \xi^{2} \xi_{x} \xi_{xx} \phi_{u}
 - 7 \xi_{t} \phi_{xxu} u
 \vspace{1mm}\\
 \quad
 - \xi^{2} \xi_{xxtt} + 4 \xi^{2} \phi_{xxxu} u
     - 7 \xi \xi_{x} \phi_{xxu} u + 4 \xi^{2} \phi \phi_{uuu} \phi_{x}
 - 6 \xi \xi_{t} \xi_{x} \xi_{xx}
+ \xi^{2} \xi_{tt} \vspace{1mm}\\
\quad
- 4 \xi \phi_{uu} \phi_{xx} u
 - 6 \xi^{2} \xi_{x} \phi_{uu} \phi_{x}
 + 6 \xi^{2} \phi_{u} \phi_{uu} \phi_{x}
  - 2 \xi^{2} \xi_{xx} \gamma u
  - 3 \xi_{xx} \phi \phi_{uu} u + 2 \xi \xi_{x} \xi_{xxx} u
  \vspace{1mm}\\
  \quad
- \xi_{t} \phi
 + 4 \xi_{x} \phi_{u} \phi_{xu} u + 4 \xi \xi_{t} \xi_{xx} \phi_{u}=0
\ea
\]

\def\CQG{Classical Quantum Grav.}
\def\IP{Inverse Problems}
\def\JPA{J. Phys. A: Math. Gen.}
\def\JPB{J. Phys. B: At. Mol. Phys.}        
\def\jpb{J. Phys. B: At. Mol. Opt. Phys.}   
\def\JPC{J. Phys. C: Solid State Phys.}     
\def\JPCM{J. Phys: Condensed Matter}        
\def\JPD{J. Phys. D: Appl. Phys.}
\def\JPE{J. Phys. E: Sci. Instrum.}
\def\JPF{J. Phys. F: Metal Phys.}
\def\JPG{J. Phys. G: Nucl. Phys.}           
\def\jpg{J. Phys. G: Nucl. Part. Phys.}     
\def\NL{Nonlinearity}
\def\PMB{Phys. Med. Biol.}
\def\RPP{Rep. Prog. Phys.}
\def\SST{Semicond. Sci. Technol.}
\def\SUST{Supercond. Sci. Technol.}
\def\AC{Acta Crystrallogr.}
\def\AM{Acta Metall.}
\def\AP{Ann. Phys., Lpz.}
\def\APNY{Ann. Phys., NY}
\def\APP{Ann. Phys., Paris}
\def\CJP{Can. J. Phys.}
\def\JAP{J. Appl. Phys.}
\def\JCP{J. Chem. Phys.}
\def\JJAP{Japan J. Appl. Phys.}
\def\JP{J. Physique}
\def\JPhCh{J. Phys. Chem.}
\def\JMMM{J. Magn. Magn. Mater.}
\def\JMP{J. Math. Phys.}
\def\JOSA{J. Opt. Soc. Am.}
\def\JPSJ{J. Phys. Soc. Japan}
\def\JQSRT{J. Quant. Spectrosc. Radiat. Transfer}
\def\NC{Nuovo Cim.}
\def\NIM{Nucl. Instrum. Methods}
\def\NP{Nucl. Phys.}
\def\PL{Phys. Lett.}
\def\PR{Phys. Rev.}
\def\PRL{Phys. Rev. Lett.}
\def\PRS{Proc. R. Soc.}
\def\PS{Phys. Scr.}
\def\PSS{Phys. Status Solidi}
\def\PTRS{Phil. Trans. R. Soc.}
\def\RMP{Rev. Mod. Phys.}
\def\RSI{Rev. Sci. Instrum.}
\def\SSC{Solid State Commun.}
\def\ZP{Z. Phys.}
\def\ARMA{Arch.\ Rat.\ Mech.\ Anal.}
\def\CMP{Commun.\ Math.\ Phys.}
\def\CPAM{Commun.\ Pure Appl.\ Math.}
\def\IP{Inverse Problems}
\def\JPA{J.\ Phys.\ A: Math.\ Gen.}
\def\NL{Nonlinearity}
\def\JMP{J.\ Math.\ Phys.}
\def\JPSJ{J.\ Phys.\ Soc.\ Japan}
\def\NC{Nuovo Cim.}
\def\PL{Phys.\ Lett.}
\def\PRA{Phys.\ Rev.\ A}
\def\PRL{Phys.\ Rev.\ Lett.}
\def\PRS{Proc.\ R.\ Soc.\ London A}
\def\PS{Phys.\ Scr.}
\def\PTRS{Phil.\ Trans.\ R. Soc.\ London A}
\def\SAM{Stud.\ Appl.\ Math.}

\label{clarkson-lp}

\end{document}